\documentclass[lettersize,journal,onecolumn]{IEEEtran}
\usepackage{amsmath,amsfonts,amsthm}
\usepackage{algorithmic}
\usepackage{algorithm}
\usepackage{array}
\usepackage[caption=false,font=normalsize,labelfont=sf,textfont=sf]{subfig}
\usepackage{textcomp}
\usepackage{stfloats}
\usepackage{url}
\usepackage{verbatim}
\usepackage{graphicx}
\usepackage{cite}
\newtheorem{assumption}{Assumption}
\newtheorem{theorem}{Theorem}
\newtheorem{lemma}{Lemma}
\newtheorem{definition}{Definition}

\newtheorem{remark}{Remark}

\newcommand{\cF}{\mathcal{F}}

\newcommand{\R}{\mathbb{R}}

\newcommand{\cS}{\mathcal{S}}

\newcommand{\order}[1]{\mathcal{O}\left(#1\right)}

\newcommand{\prt}[1]{\left(#1\right)}
\newcommand{\brk}[1]{\left[#1\right]}
\newcommand{\crk}[1]{\left\{#1\right\}}

\newcommand{\orderi}[1]{\mathcal{O}(#1)}
\newcommand{\prti}[1]{(#1)}
\newcommand{\brki}[1]{[#1]}
\newcommand{\crki}[1]{\{#1\}}

\newcommand{\inpro}[1]{\left\langle #1 \right\rangle}

\newcommand{\condE}[2]{\E\brk{#1\middle|#2}}
\newcommand{\condEi}[2]{\E[#1|#2]}
\newcommand{\norm}[1]{\left\Vert #1 \right\Vert}
\newcommand{\normi}[1]{\Vert #1 \Vert}
\newcommand{\bs}[1]{\boldsymbol{#1}}

\newcommand{\E}{\mathbb{E}}

\newcommand{\sumn}{\sum_{i=1}^n}

\newcommand{\sigfn}{\sigma_f^\star}

\usepackage{booktabs}
\usepackage{tabularx}
\usepackage[flushleft]{threeparttable}
\usepackage{makecell} 
\usepackage{pifont}
\usepackage{float}
\usepackage{multirow}

\usepackage{url}            			
\usepackage{authblk}

\makeatletter
\def\@buildmath#1{%
    \expandafter\def\csname bb#1\endcsname{\mathbb{#1}}%
    \expandafter\def\csname bs#1\endcsname{\boldsymbol{#1}}%
    \expandafter\def\csname cal#1\endcsname{\mathcal{#1}}%
    \expandafter\def\csname frak#1\endcsname{\mathfrak{#1}}%
    \expandafter\def\csname rm#1\endcsname{\mathrm{#1}}%
    \expandafter\def\csname bf#1\endcsname{\mathbf{#1}}%
}
\def\@buildmathletters#1{%
    \ifx#1\relax
    \else
        \@buildmath{#1}%
        \expandafter\@buildmathletters
    \fi
}
\@buildmathletters ABCDEFGHIJKLMNOPQRSTUVWXYZabcdefghijklmnopqrstuvwxyz\relax
\makeatother

\usepackage{color}

\definecolor{cuhkpl}{RGB}{152,24,147}

\definecolor{beanGreen}{RGB}{255, 1, 200}

\usepackage[colorlinks=true, linkcolor=blue, citecolor=blue]{hyperref}

\begin{document}

\title{On Federated Averaging with Partial Participation under a General Variance Condition}

\author{Ziming Wang, 
Kun Huang,
Karl Henrik Johansson,
and Shi Pu
\thanks{Ziming Wang and Kun Huang contribute equally to this work. Ziming Wang and Shi Pu are with the School of Data Science, The Chinese University of Hong Kong, Shenzhen, China. Kun Huang and Karl Henrik Johansson are with the Department of Decision and Control Systems, School of Electrical Engineering and Computer Science, KTH Royal Institute of Technology, and Digital Futures, Stockholm, Sweden. \tt\small (emails: zimingwang2@link.cuhk.edu.cn, kunhuang@kth.se, kallej@kth.se, pushi@cuhk.edu.cn)}
}



\maketitle

\begin{abstract}
In this paper, we study the convergence of Federated Averaging (FedAvg) with partial participation under a general variance condition. In each communication round, the server samples $S$ of the $n$ clients uniformly at random without replacement. The stochastic gradients are assumed to satisfy the so-called ABC condition, which allows their variance to grow with the local objective gaps. For smooth nonconvex objectives, we establish a finite-time stationarity bound in terms of the heterogeneity measure $\sigma_f^\star:=f^\star-\frac{1}{n}\sum_{i=1}^n f_i^\star$ and the client sampling factor $I_S:=(n-S)/(n-1)$. Here, $f^\star$ and $f_i^\star$ denote the infima of the global average objective and client $i$'s local objective, respectively. The bound exhibits a linear speedup in the number of sampled clients $S$ and recovers the full-participation result when $S=n$.
Under diminishing stepsizes, we further prove that the gradient norm converges to zero both in expectation and almost surely. Notably, compared with existing analyses of FedAvg with partial participation, our results require no additional data heterogeneity assumption such as bounded gradient dissimilarity (BGD), and the sampling-induced term in our bound is not amplified by the number of local updates at the leading order. The importance of removing the BGD assumption is illustrated by an example that is smooth
and lower bounded but violates BGD, yet is covered by our analysis.
\end{abstract}


\section{Introduction}

Federated learning (FL) trains a shared model across clients while keeping each client's raw data local. It has been applied to autonomous driving~\cite{chen2021bdfl,zeng2022federated}, large-scale training~\cite{sani2025photon}, and natural language processing~\cite{yang2018applied,xu2023federated}. In FL, the clients collaboratively solve
\begin{equation}
	\label{eq:P}
	\min_{x \in \mathbb{R}^d} f(x) := \frac{1}{n} \sum_{i=1}^n f_i(x),
\end{equation}
where $n$ is the number of clients and \(f_i:\R^d\to\R\) is the local objective available only to client $i\in[n]:=\crki{1,\ldots,n}$. Typical examples of $f_i$ include the expected risk $f_i(x)=\mathbb{E}_{\xi\sim\mathcal{D}_i}[F_i(x;\xi)]$ and the empirical risk $f_i(x)=\frac{1}{m_i}\sum_{j=1}^{m_i} f_{ij}(x)$, where $\mathcal{D}_i$ is the local data distribution, $F_i(x;\xi)$ is the loss of model $x$ on sample $\xi$, and $f_{ij}(x):=F_i(x;\xi_{ij})$ is the loss on the $j$-th local sample $\xi_{ij}$.
Since client data typically follow different distributions, the local objectives $f_i$ generally differ across clients. This data heterogeneity is a central concern in FL.

A classical algorithm for solving~\eqref{eq:P} is Federated Averaging (FedAvg)~\cite{mcmahan2017communication}. We consider FedAvg with partial client participation: in each communication round, the server samples $S$ of the $n$ clients uniformly at random without replacement. Each sampled client performs $Q$ local stochastic gradient updates starting from the current global model, and the server averages the resulting local models. 
Partial participation is motivated by limited server bandwidth and
intermittent client availability~\cite{bonawitz2019towards,kairouz2021advances}.

Convergence analyses for FedAvg typically impose conditions on both the stochastic gradient noise and the data heterogeneity. Let $g_i(x;\xi_i)$ be an unbiased stochastic gradient of $f_i(x)$. A commonly used noise condition is the relaxed growth condition~\cite{bottou2018optimization},
\begin{equation}
    \label{eq:bv_intro}
		\E_{\xi_i}\brk{\norm{g_i(x;\xi_i) - \nabla f_i(x)}^2}\leq \sigma^2 + \eta^2\norm{\nabla f_i(x)}^2,\; \forall x\in\R^d,
\end{equation}
for some constants $\sigma,\eta \ge 0$.
In particular, its special case $\eta=0$ (bounded variance) has been widely adopted~\cite{yang2021achieving,karimireddy2020scaffold}.

To model and control data heterogeneity, one typically imposes the bounded gradient dissimilarity (BGD) condition (see, e.g., \cite{woodworth2020minibatch,zhang2021fedpd}):
\begin{equation}
    \label{Zb}
    \frac{1}{n}\sum_{i=1}^n
    \|\nabla f_i(x)-\nabla f(x)\|^2
    \le
    \zeta^2+\psi^2\|\nabla f(x)\|^2,
    \quad \forall x\in\mathbb R^d,
\end{equation}
where $\zeta,\psi\ge0$ are constants. The left-hand side of~\eqref{Zb} vanishes for all $x$ if and only if all $f_i$ have identical gradients. The constants $\zeta$
and $\psi$ quantify the level of data heterogeneity. 

Both conditions can be restrictive. The relaxed growth condition~\eqref{eq:bv_intro} can fail {for} simple finite-sum problems~\cite{huang2024distributed,khaled2020better}. Moreover, BGD can fail even when all $f_i$ are smooth and lower bounded: Theorem~\ref{thm:BGD} constructs such local objectives for which no finite pair $(\zeta,\psi)$ satisfies~\eqref{Zb}.  
To our knowledge, however, existing analyses of FedAvg with partial participation still rely on BGD (or a stronger condition) together with bounded variance. Specifically, client sampling introduces an additional error that scales with the gradient dissimilarity, and these analyses (e.g.,~\cite{yang2021achieving}) control this error through BGD. Beyond FedAvg, the methods in~\cite{huang2024stochastic,you2025fedsum} avoid BGD under partial participation, but they still assume bounded variance. More importantly, they require clients to retain auxiliary variables across rounds. In contrast, FedAvg is stateless and well suited to cross-device FL, where clients participate infrequently and stored optimization variables can become stale~\cite{kairouz2021advances}.

The above observations motivate us to establish convergence guarantees for FedAvg with partial participation under a more general noise model and without BGD-type assumptions.
The ABC condition~\cite{khaled2020better,lei2019stochastic,huang2024distributed} provides such a noise model:
\begin{equation}
	\label{eq:abc_intro}
	\E_{\xi_i}\brk{\norm{g_i(x;\xi_i) - \nabla f_i(x)}^2}
    \leq C\brk{f_i(x)-f_i^\star}+D,
\end{equation}
where $C,D\ge0$ are constants and $f_i^\star := \inf_{x} f_i(x)$; see Assumption~\ref{A2}. For smooth and lower bounded $f_i$, the squared gradient norm is controlled by the local objective gap, so \eqref{eq:bv_intro} implies \eqref{eq:abc_intro}. The converse need not hold~\cite{huang2024distributed,khaled2020better}. Thus, the ABC condition covers a strictly broader class of stochastic gradient oracles.


Regarding data heterogeneity, BGD can be avoided. Let
$f^\star:=\inf_{x}f(x)$ and $\sigfn:=f^\star-\frac{1}{n}\sum_{i=1}^n f_i^\star$. If each $f_i$ is $L$-smooth and lower bounded, the left-hand side of~\eqref{Zb} is at most $2L[f(x)-f^\star + \sigfn]$; see Lemma~\ref{AL1}. Since $\sigfn$ is finite whenever each $f_i$ is lower bounded, this bound remains valid even when BGD fails~\cite{huang2024distributed}.

In this paper, we prove that FedAvg with partial participation converges to stationarity under the ABC condition~\eqref{eq:abc_intro}, without BGD-type assumptions.
The main difficulty is that the aggregated update is affected by two sources of randomness. The first is the stochastic gradient noise along the local trajectories, which scales with the local objective gaps. The second is client sampling, whose error scales with the gradient dissimilarity. 
Without BGD, this dissimilarity is bounded in terms of $2L[f(x)-f^\star + \sigfn]$, which may grow even when $\normi{\nabla f(x)}^2$ stays bounded. 
We separate the two sources by assigning a virtual local trajectory to every client, and bound the gradient dissimilarity by $f(x)-f^\star$ and $\sigfn$. The resulting bound quantifies the effect of partial participation through the sampling factor $I_S:=(n-S)/(n-1)$, which vanishes under full participation. Moreover, the client-sampling error in our bound is not amplified by the number of local updates.

\begin{table*}[htbp]
		\centering
		\begin{threeparttable}
			\caption[Comparison of convergence results of FedAvg for smooth nonconvex objectives]{Comparison of convergence results of FedAvg for smooth nonconvex objectives.}
			\label{tab:convergence_comparison}
			\scriptsize 
			\setlength{\tabcolsep}{2pt} 
				\begin{tabular*}{0.85\textwidth}{@{\extracolsep{\fill}}l c c c c c@{}}
				\toprule
				\multirow{2}{*}{Paper} & \multirow{2}{*}{Noise} & \multicolumn{2}{c}{Heterogeneity} & \multirow{2}{*}{Rate} & \multirow{2}{*}{a.s.} \\ \cmidrule(lr){3-4}
				& & Assumption & Measure & & \\ \midrule
				\multicolumn{6}{l}{\textbf{Full Participation ($\boldsymbol{S=n}$)}} \\
				\cite{yu2019linear} & $(\sigma^2, 0)$ & $(\zeta^2, 0)$-BGD & $\zeta^2$ & $\order{\frac{1+\sigma^2}{\sqrt{nQT}}} $ & N \\
				
				\cite{yang2021achieving} & $(\sigma^2, 0)$ & $(\zeta^2, 0)$-BGD & $\zeta^2$ & $\order{\frac{1+\sigma^2}{\sqrt{nQT}}}$ & N \\
				
				\cite{wang2020tackling} & $(\sigma^2, 0)$ & $(\zeta^2, \psi^2)$-BGD & $\zeta^2,\psi^2$ & $\order{\frac{1+\sigma^2}{\sqrt{nQT}}}$ & N \\
				
				\cite{Qin03092025} & $(\sigma^2, \eta^2)$ & $(\zeta^2, \psi^2)$-BGD & $\zeta^2,\psi^2$ & $\order{\frac{1+\sigma^2}{\sqrt{nQT}}}$ & N \\
                
				\cite{huang2024distributed} & ABC & None & ${\sigfn}^{(a)}$ & $\order{\sqrt{\frac{D+C\sigfn+C}{nQT}}}$ & {Y} \\ 
                \midrule
				
				\multicolumn{6}{l}{\textbf{Partial Participation ($\boldsymbol{S< n}$)}} \\

				 \cite{yang2021achieving} & $(\sigma^2, 0)$ & $(\zeta^2, 0)$-BGD & $\zeta^2$ & $\order{\frac{1+\sigma^2}{\sqrt{SQT}} + \frac{I_S \zeta^2\sqrt{Q}}{\sqrt{ST}}}^{(b)}$
                  & N\\


				\makecell[c]{\textbf{This} \textbf{Work}} & \textbf{ABC} & \textbf{None} & \textbf{$\sigfn$} & \makecell[c]{$\bs{\order{\sqrt{\frac{D+C\sigfn+C}{SQT}} + \sqrt{\frac{(1+\sigfn)I_S }{ST}}}}$} & \textbf{Y} \\
				\bottomrule
			\end{tabular*}
            \begin{tablenotes}
            \item Noise: $(\sigma^2,\eta^2)$ means $\E\normi{g_i(x;\xi)-\nabla f_i(x)}^2\le\sigma^2+\eta^2\normi{\nabla f_i(x)}^2$; ABC refers to~\eqref{eq:abc_intro}. Heterogeneity: $(\zeta^2,\psi^2)$-BGD refers to~\eqref{Zb}. Rate: only leading-order terms in $T$ are reported. Problem-dependent constants and numerical constants are omitted. a.s.: almost sure convergence under $\sum_t\alpha_t=\infty$ and $\sum_t\alpha_t^2<\infty$.
				\item $^{(a)}$ $\sigfn = f^\star - \frac{1}{n}\sum f_i^\star$ with $f^\star = \inf_x f(x)$ and $f_i^\star = \inf_x f_i(x)$.
                \item $^{(b)}$ $I_S =(n-S)/(n-1)$ is the sampling factor.
			\end{tablenotes}
		\end{threeparttable}
	\end{table*}

\subsection{Related Works}


Most convergence analyses of FedAvg control data heterogeneity through the BGD condition~\eqref{Zb}; see, e.g.,~\cite{woodworth2020minibatch,yu2019linear,wang2020tackling,yang2021achieving,Qin03092025}.
Other works avoid BGD by adopting alternative heterogeneity measures, such as the average drift at the optimum~\cite{wang2022unreasonable} and $\sigfn$~\cite{li2020fedavg,huang2024distributed}. The measure
$\sigfn$ is used in~\cite{li2020fedavg} under strong convexity and bounded second moments of the stochastic gradients. It is shown in~\cite{huang2024distributed} that $\sigfn$ suffices to
analyze nonconvex FedAvg under the ABC condition, but only for full
participation.
A different line of work modifies FedAvg to mitigate heterogeneity: SCAFFOLD~\cite{karimireddy2020scaffold} corrects client drift with control variates, FedProx~\cite{li2020fedprox} adds a proximal term to the local objectives, and FedPD~\cite{zhang2021fedpd} adopts a primal-dual formulation. 

The effect of client sampling on convergence has been widely studied. For strongly convex objectives, FedAvg has been analyzed under client sampling both with and without replacement~\cite{li2020fedavg}. For
smooth nonconvex objectives, a rate of $\orderi{\sqrt{Q/(ST)}}$ has been
established for FedAvg with partial participation~\cite{yang2021achieving}. Other works design optimal
client sampling rules that better approximate the full-client update~\cite{chen2022optimal} and handle arbitrary participation patterns~\cite{you2025fedsum}. These results are mainly developed under bounded variance or BGD-type assumptions, whereas our analysis requires neither.


FL is closely related to decentralized optimization, where clients exchange information with their neighbors over a connected network instead of through a central server.
Classical decentralized methods include DSGD~\cite{lian2017can,pu2021sharp},
gradient tracking~\cite{xu2015augmented,di2016next,nedic2017achieving,pu2021distributed},
and exact diffusion/$D^2$~\cite{yuan2018exact,tang2018d,huang2022improving}.
Their local-update variants are studied in~\cite{koloskova2020unified,huang2023computation,liu2024decentralized,alghunaim2024local}, and local momentum tracking further attains the optimal communication complexity~\cite{huang2025distributed}.
The stochastic versions of these
methods are typically analyzed under bounded variance or the relaxed
growth condition~\eqref{eq:bv_intro}. Decentralized methods under the ABC condition~\eqref{eq:abc_intro} have been studied in~\cite{huang2024nmap,huang2026accelerated}, but without local
updates. 

\subsection{Main Contribution}

The main contribution of this work is three-fold.

First, we establish a finite-time stationarity guarantee for FedAvg with partial participation under the ABC condition~\eqref{eq:abc_intro} without BGD-type assumptions. In particular, with $S$ sampled clients per round, $Q$ local updates, and
$T$ total communication rounds, the leading terms of the bound on $\min_{0\le t<T}\E\normi{\nabla f(x_t)}^2$ are
\small
\begin{equation*}
	\order{\sqrt{\frac{L\Delta_f(D+C\sigfn+C\Delta_f)}{SQT}}+L\prt{\Delta_f+\sqrt{\Delta_f\sigfn}}\sqrt{\frac{I_S}{ST}}},
\end{equation*}\normalsize
where $\Delta_f:=f(x_0)-f^\star$. The first term stems from the stochastic gradient noise and exhibits a linear speedup in $S$ and $Q$.
The second term captures the effect of client sampling and is independent of $Q$. By comparison, the corresponding sampling term in~\cite{yang2021achieving} grows as $\sqrt{Q}$ under their stepsize choice. Under full participation, $I_S=0$ and the bound recovers the $\orderi{1/\sqrt{nQT}}$ rate
of~\cite{huang2024distributed}. To the best of our knowledge, this is the first guarantee for FedAvg with partial participation that
requires neither bounded variance nor BGD-type assumptions; see
Table~\ref{tab:convergence_comparison}.

Second, under diminishing stepsizes, we show that $\normi{\nabla f(x_t)}$ converges to $0$ both in expectation and almost surely. To the best of our knowledge, this is the first almost sure convergence result for FedAvg with partial participation in the nonconvex setting. It establishes an asymptotic stationarity guarantee for the last iterate and complements the finite-time bound above.

Third, we show that the BGD condition~\eqref{Zb} is indeed restrictive: it
excludes certain problems that satisfy the smoothness and lower
boundedness assumptions. Specifically, Theorem~\ref{thm:BGD} constructs a one-dimensional example with two clients whose local objectives are smooth and lower bounded. However, no finite pair $(\zeta,\psi)$ satisfies~\eqref{Zb} in this example, even though the data heterogeneity measured by $\sigfn$ is at most $1/2$. Hence, within this setting, our analysis covers a strictly larger class of problems compared to BGD-based convergence guarantees.

\section{FedAvg with Partial Participation: Algorithm and Preliminary Analysis}\label{sec:algorithm}

In this section, we describe FedAvg with partial participation, state the standing assumptions, introduce the heterogeneity measure, and
develop the auxiliary estimates used in its convergence analysis.

\subsection{Algorithm and Virtual Trajectories}\label{subsec:algorithm}
Algorithm~\ref{Al1} summarizes FedAvg with partial participation. In each communication round $t$, the server samples a subset $\mathcal S_t$ of $S$ clients uniformly at random without replacement and sends the current global model $x_t$ to them. Each sampled client starts from $x_t$, performs $Q$ local stochastic gradient steps with the common stepsize $\alpha_t$, and returns its final local model. The server then averages these models to obtain $x_{t+1}$. The unselected clients remain idle: they neither perform local computation nor keep any local state across rounds, since a client selected in a later round $k>t$ simply restarts from $x_k$. This participation model reflects practical constraints such as intermittent connectivity and limited client-side resources~\cite{mcmahan2017communication,luo2024adaptive}.
\begin{algorithm}[h]
	\caption{FedAvg with Partial Participation}
	\label{Al1}
	\begin{algorithmic}[1]
		\STATE \textbf{Server input:} initial $x_0$, number of communication rounds $T$ and number of participating clients $S$.
		\STATE \textbf{Client input:} number of local updates $Q$ and stepsize sequence $\crki{\alpha_t}_{t\geq 0}$.
		\FOR{$t = 0,1,\dots,T-1$}
		\STATE Server samples $\mathcal{S}_t$ uniformly without replacement, $|\mathcal{S}_t| = S$
		\STATE Send $x_{t}$ to all  clients $i\in \mathcal{S}_t$
		\FOR{Client $i \in \mathcal{S}_t$ \textbf{in parallel}}
		\STATE Receive $x_t$ from server and set $x_{i,t}^0 = x_t$.
		\FOR{$\ell = 0, 1, \dots, Q - 1$}
		\STATE Query a stochastic gradient $g_i(x_{i,t}^\ell; \xi_{i,t}^\ell)$.
		\STATE Update $x_{i,t}^{\ell+1} = x_{i,t}^\ell - \alpha_t g_i(x_{i,t}^\ell; \xi_{i,t}^\ell)$.
		\ENDFOR
		\ENDFOR
		\STATE Server aggregates $x_{t+1} = \frac{1}{S} \sum_{i\in \mathcal{S}_t} x_{i,t}^Q$.
       \ENDFOR
       \STATE Output $x_T$.
	\end{algorithmic}
\end{algorithm}

For each selected client \(i\in\mathcal S_t\), its local recursion can be written as
\begin{equation}
	x_{i,t}^{\ell+1} = x_{i,t}^\ell - \alpha_t g_i(x_{i,t}^\ell; \xi_{i,t}^\ell) = x_t - \alpha_t \sum_{j=0}^\ell g_i(x_{i,t}^j; \xi_{i,t}^j),
\end{equation}
where $ \ell \in \{0, 1, \dots, Q-1\}$, $t \ge 0$ and $\xi_{i,t}^{\ell}$ denotes the random sample\footnote{For instance, a data sample or mini-batch drawn from $\mathcal D_i$.} used by the stochastic gradient oracle of client~$i$ at local step~$\ell$ of round~$t$.
Consequently, the server update can be expressed as
\begin{equation}
    x_{t+1}=x_t-\alpha_t G_t,\quad G_t:=\frac{1}{S}\sum_{i\in \mathcal{S}_t}\sum_{\ell=0}^{Q-1}
g_i(x_{i,t}^{\ell};\xi_{i,t}^{\ell}).
\end{equation}

For the analysis, we associate each client with a virtual local trajectory at every communication round, regardless of whether the client is selected. Similar ideas also appear in~\cite{li2020fedavg,yang2021achieving}. To this end, we regard the samples $\xi_{i,t}^{\ell}$ as defined for all clients $i\in[n]$, including those not selected at round~$t$. However, only the samples of the selected clients enter the algorithm.
For each client $i\in[n]$, define the virtual local trajectory
\[
    \widetilde x_{i,t}^{0}=x_t,
    \quad
    \widetilde x_{i,t}^{\ell+1}
    =\widetilde x_{i,t}^{\ell}-\alpha_t\, g_i(\widetilde x_{i,t}^{\ell};\xi_{i,t}^{\ell}),
    \quad \ell=0,\ldots,Q-1.
\]
By construction, the virtual trajectories extend the actual local iterates to all clients:
\[
    \widetilde x_{i,t}^{\ell+1}=
    \begin{cases}
        x_{i,t}^{\ell+1}, &  i\in\mathcal S_t,\\
        \widetilde x_{i,t}^{\ell}-\alpha_t\, g_i(\widetilde x_{i,t}^{\ell};\xi_{i,t}^{\ell}), &  i\notin\mathcal S_t.
    \end{cases}
\]
Consequently, the aggregated direction can be
written as
\[
    G_t=\frac{1}{S}\sum_{i\in\mathcal S_t}\sum_{\ell=0}^{Q-1} g_i(\widetilde x_{i,t}^{\ell};\xi_{i,t}^{\ell}),
\]
that is, $G_t$ is the average of $S$ directions sampled without replacement from the finite population
$\crki{\sum_{\ell=0}^{Q-1} g_i(\widetilde x_{i,t}^{\ell};\xi_{i,t}^{\ell})}_{i=1}^n$ of potential
client directions. This representation is the basis of the sampling analysis in Section~\ref{subsec:prelim}.

\subsection{Assumptions}\label{subsec:assumption}

	We impose three standing assumptions: smoothness and lower boundedness of the local objectives (Assumption~\ref{A1}), unbiased stochastic gradients satisfying the ABC condition (Assumption~\ref{A2}), and uniform client sampling without replacement (Assumption~\ref{A3}). The ABC condition~\cite{khaled2020better,lei2019stochastic} covers a strictly larger class of stochastic gradients than the bounded variance and relaxed growth conditions~\cite{huang2024distributed}. Unlike existing analyses of FedAvg with partial participation, we impose no additional data heterogeneity assumption such as the BGD condition~\eqref{Zb}, which may fail under Assumption~\ref{A1} (Theorem~\ref{thm:BGD}).
	\begin{assumption}
		\label{A1}
		Each $f_i$ is $L$-smooth and lower bounded, i.e.,
		\begin{equation}
			\begin{aligned}
				&\|\nabla f_i(x) - \nabla f_i(y)\| \le L \|x-y\|,\\ 
				&f_i(x) \geq f_i^\star :=\inf_{z\in\R^d} f_i(z)> -\infty, \quad \forall x,y\in\R^d.
			\end{aligned}
		\end{equation}
	\end{assumption}

Assumptions~\ref{A2} and~\ref{A3} are stated with respect to the following filtrations.

\begin{definition}
\label{def:filtration}
Let $\mathcal F_t :=\sigma (x_0,\,\{\mathcal S_s\}_{0\le s<t},\, \{\xi_{i,s}^{\ell}: i\in[n],0\le \ell<Q,0\le s<t\})$
be the information available at the beginning of round $t$. In particular, $x_t$ is $\mathcal F_t$-measurable.
Within round $t$, define $\mathcal F_t^0:=\mathcal F_t$ and $\mathcal{F}_t^{\ell+1}:= \sigma(\mathcal F_t^{\ell},\{\xi_{i,t}^{\ell}\}_{i=1}^n)$ for $\ell=0,\ldots,Q-1$.
The sampling set is revealed only after all round-\(t\) oracle variables, i.e., $\mathcal F_{t+1}=\sigma(\mathcal F_t^Q,\mathcal S_t)$.
Thus,
\[
    \mathcal F_t
    =
    \mathcal F_t^{0}
    \subseteq
    \mathcal F_t^{1}
    \subseteq\cdots\subseteq
    \mathcal F_t^{Q}
    \subseteq
    \mathcal F_{t+1}.
\]
Moreover, we use
$    \mathbb E_t[\cdot]
    :=
    \mathbb E[\cdot\mid\mathcal F_t]$.
\end{definition}

Since \(\widetilde x_{i,t}^0=x_t\) and the oracle $g_i$
is a measurable function of its arguments, every virtual iterate
\(\widetilde x_{i,t}^{\ell}\) is \(\mathcal F_t^{\ell}\)-measurable, and every
\(g_i(\widetilde x_{i,t}^{\ell};\xi_{i,t}^{\ell})\) is \(\mathcal F_t^{\ell+1}\)-measurable.
In particular, the round-\(t\) virtual trajectories $\crki{\widetilde x_{i,t}^{\ell}: i\in[n],\,0\le\ell\le Q}$ are \(\mathcal F_t^{Q}\)-measurable: they are determined by $x_t$ and the round-$t$ samples, and do not involve the round-$t$ sampling set $\mathcal S_t$. This is the purpose of the ordering in Definition~\ref{def:filtration}: $\mathcal F_t^Q$ contains all round-$t$ samples but not $\mathcal S_t$. The ordering is used only for the analysis and does not contradict Algorithm~\ref{Al1}, where $\mathcal S_t$ is drawn before the local updates. This is because the virtual trajectories do not involve $\mathcal S_t$ and $\mathcal S_t$ is independent of the samples (Assumption~\ref{A3}), revealing $\mathcal S_t$ after the samples does not change the joint distribution of $\mathcal S_t$ and the samples.

\begin{assumption}
\label{A2}
Each client $i$ has access to an unbiased stochastic gradient of $f_i$ whose variance satisfies the ABC condition along the virtual trajectories. Formally, there exist constants $C,D\ge0$ such that, for every $t\ge0$, $i\in[n]$, and $\ell\in\{0,\ldots,Q-1\}$, the stochastic gradient $g_i(\widetilde{x}_{i,t}^{\ell};\xi_{i,t}^{\ell})$
satisfies $\condEi{g_i(\widetilde{x}_{i,t}^{\ell};\xi_{i,t}^{\ell})}{\mathcal F_t^{\ell}}=\nabla f_i(\widetilde{x}_{i,t}^{\ell})$  and
\begin{equation}\label{ABC}
   \condE{\norm{g_i(\widetilde{x}_{i,t}^{\ell};\xi_{i,t}^{\ell})-\nabla f_i(\widetilde{x}_{i,t}^{\ell})}^2}{\mathcal F_t^{\ell}}
    \leq
    C\brk{f_i(\widetilde x_{i,t}^{\ell})-f_i^\star}+D
\end{equation}
almost surely (a.s.). Moreover, conditionally on $\mathcal F_t^{\ell}$, the noises $\{\xi_{i,t}^{\ell}\}_{i=1}^n$ are independent across clients.
\end{assumption}
\begin{remark}
    At the beginning of each round, all virtual local trajectories start from the same point, $\widetilde x_{i,t}^{0}=x_t$. Averaging~\eqref{ABC} over clients therefore gives,
    \begin{align*}
    &\frac{1}{n}\sum_{i=1}^n
    \E_t\brk{
        \norm{g_i(x_t;\xi_{i,t}^{0})-\nabla f_i(x_t)}^2
    }\\
    &\leq
    C[f(x_t)-f^\star]+C\sigfn+D.
    \end{align*}
    This upper bound depends on the global objective gap through $C[f(x_t)-f^\star]$ and on data heterogeneity through $C\sigfn$.
    The bounded-variance assumption is the special case $C=0$ and $D=\sigma^2$, for which the bound contains no explicit heterogeneity term.
\end{remark}

\begin{assumption}
\label{A3}
Conditionally on \(\mathcal F_t^Q\), the server samples
\(\mathcal S_t\subseteq[n]\) uniformly at random without replacement, i.e.,
\begin{align*}
    \mathbb P(\mathcal S_t=\mathcal A\mid\mathcal F_t^Q)= \frac{1}{\binom{n}{S}},\quad \forall \mathcal A\subseteq[n],\quad |\mathcal A|=S.
\end{align*}
\end{assumption}
Let \(\mathbb{I}_{i,t}:=\mathbf 1_{\{i\in\mathcal S_t\}}\). By Assumption~\ref{A3},
\[
    \condE{\mathbb{I}_{i,t}}{\mathcal F_t^Q} = \frac{S}{n},
    \quad
    \condE{\mathbb{I}_{i,t}\cdot\mathbb{I}_{j,t}}{\mathcal F_t^Q}= \frac{S(S-1)}{n(n-1)}
    \quad i\neq j.
\]
Consequently, by the tower property, for any integrable \(\mathcal F_t^Q\)-measurable random variables \(Z\) and \(\{a_i\}_{i=1}^n\),
\begin{equation}\label{eq:indep}
	\begin{aligned}
		 \E_t\brk{\mathbb I_{i,t}Z}&=\frac{S}{n}\mathbb E_t[Z],\\
		\E_t\brk{\mathbb I_{i,t}\cdot\mathbb I_{j,t}Z}
		&=\frac{S(S-1)}{n(n-1)}\mathbb E_t\brk{Z}, \quad i\neq j,\\
		\condE{\frac{1}{S}\sum_{i\in\mathcal S_t}a_i}{\mathcal{F}_t^Q} &=\frac{1}{n}\sum_{i=1}^n a_i .
	\end{aligned}
\end{equation}
Since the round-\(t\) virtual trajectories are \(\mathcal F_t^{Q}\)-measurable, the randomness of round-\(t\) client sampling can be separated from that of the oracles through \eqref{eq:indep}, while the oracle noise is handled through the filtration \(\{\mathcal F_t^{\ell}\}_{\ell=0}^{Q}\).
We also define the sampling factor
\begin{equation}\label{eq:IS}
    I_S:=\frac{n-S}{n-1},\quad S\in\{1,\dots,n\},\quad n\geq 2.
\end{equation}

\subsection{Heterogeneity Measure and a Counterexample to the BGD Condition}\label{subsec:hetero}
    Under Assumption~\ref{A1}, the global objective $f=\frac{1}{n}\sum_{i=1}^n f_i$ is also $L$-smooth and lower bounded, since $f(x)\ge\frac{1}{n}\sum_{i=1}^n f_i^\star$ for all $x\in\R^d$. Let $f^\star:=\inf_{x\in\R^d}f(x)$, and define the heterogeneity measure
    \begin{equation}\label{def:sigfn}
        \sigfn:=f^\star-\frac{1}{n}\sum_{i=1}^n f_i^\star.
    \end{equation}
    Equivalently, $\sigma_f^\star=\inf_{x\in\R^d}\frac{1}{n}\sum_{i=1}^n[f_i(x)-f_i^\star]$ is the infimum of the average local objective gap. Hence it is nonnegative and vanishes when the local objectives share a common minimizer. The following lemma, taken from~\cite[Lemma 1]{huang2024distributed}, shows that $\sigma_f^\star$, together with the functional gap $f(x)-f^\star$, controls the gradient dissimilarity and avoids the need for additional data heterogeneity assumptions.

    \begin{lemma}\label{AL1}
        Let Assumption~\ref{A1} hold. Then, for every $x \in \R^d$
    \begin{equation}
        \frac{1}{n} \sum_{i=1}^n \norm{\nabla f_i(x) - \nabla f(x)}^2 \leq 2L \brk{ f(x) - f^\star } + 2L \sigma_f^\star.
    \end{equation}
    \end{lemma}

    Compared with the BGD condition~\eqref{Zb}, $2L\sigma_f^\star$ plays the role of $\zeta^2$, and the term $2L[f(x)-f^\star]$ replaces $\psi^2\normi{\nabla f(x)}^2$. Unlike the constants $\zeta$ and $\psi$, whose existence must be assumed, $\sigma_f^\star$ is always finite under Assumption~\ref{A1}. The following theorem constructs a problem that satisfies Assumption~\ref{A1} but violates~\eqref{Zb} for every pair $(\zeta,\psi)$. Hence, Lemma~\ref{AL1} remains applicable where the BGD condition fails.

    \begin{theorem}\label{thm:BGD}
    Let $n=2$, $d=1$, $r(x):=\sqrt{1+x^2}$, and
    \begin{align}\notag
    	f_1(x)&=r(x)+\frac{1}{2}r(x)\sin\sqrt{r(x)},\\
    	f_2(x)&=r(x)-\frac{1}{2}r(x)\sin\sqrt{r(x)}.\notag
    \end{align}
    Then $f_1$ and $f_2$ satisfy Assumption~\ref{A1} with $L=9/4$, and $\sigfn\le 1/2$. However, for any constants $\zeta,\psi\ge0$, there exists $x\in\R$ at which the BGD condition~\eqref{Zb} is violated.
    \end{theorem}

    \begin{proof}
    Let $u(x):=\frac{1}{2}r(x)\sin\sqrt{r(x)}$, so that $f_1=r+u$, $f_2=r-u$, and $f=\frac{1}{2}(f_1+f_2)=r$. Note that $r(x)\ge1$ for all $x\in\R$. We have 
    \begin{align*}
        r'(x) = \frac{x}{r(x)} \in(-1,1), \quad r''(x) = \frac{1}{r(x)^3} \in(0,1].
    \end{align*}

    \textbf{Lower boundedness.} Since $|u|\le r/2$, we have $f_i\ge r/2\ge 1/2$ for $i=1,2$, and hence $f_i^\star\ge 1/2$. Moreover, $f^\star=\inf_{x\in\R}r(x)=1$, so 
    \begin{align*}
        \sigfn = f^\star-\frac{1}{2}(f_1^\star+f_2^\star)\leq \frac{1}{2}.
    \end{align*}

    \textbf{Smoothness.} Let $U(R):=\frac{R}{2}\sin\sqrt{R}$ for $R\ge1$, so that $u=U\circ r$. Direct computation gives
    \begin{align*}
    	U'(R)&=\frac{1}{2}\sin\sqrt{R}+\frac{\sqrt{R}}{4}\cos\sqrt{R},\\
    	U''(R)&=\frac{3}{8\sqrt{R}}\cos\sqrt{R}-\frac{1}{8}\sin\sqrt{R}.
    \end{align*}
    Hence $|U'(R)|\le 1/2+\sqrt{R}/4$ and $|U''(R)|\le 1/2$ for $R\ge1$. Since $u' = U'(r) r'$, $u''=U''(r)(r')^2+U'(r)r''$, $|r'|<1$, $r''=r^{-3}$.
    We obtain
    \[
        |u''(x)|\le\frac{1}{2}+\frac{1}{r^3}\prt{\frac{1}{2}+\frac{\sqrt{r}}{4}}\le\frac{1}{2}+\frac{1}{2}+\frac{1}{4}=\frac{5}{4},\qquad r\geq 1.
    \]
    Therefore 
    \begin{align*}
        |f_i''(x)| \leq |r''(x)| + |u''(x)| \leq \frac{1}{r^3} + \frac{5}{4} \leq \frac{9}{4},\qquad i=1,2.
    \end{align*}
    Hence, $f_1$ and $f_2$ are $L$-smooth with $L=9/4$. Together with the lower boundedness, Assumption~\ref{A1} holds.

    \textbf{Failure of BGD.} Since $\nabla f_1(x) - \nabla f(x) = u'(x)$, $\nabla f_2(x) - \nabla f(x) = -u'(x)$, $\nabla f(x) = r'(x)$,
    we have 
    \begin{align*}
        \frac{1}{2}\sum_{i=1}^2 \norm{\nabla f_i(x) - \nabla f(x)}^2 &= \left|u'(x)\right|^2,\\
         \zeta^2 + \psi^2 \left|r'(x)\right|^2 &\le \zeta^2 + \psi^2.
    \end{align*}

    We now construct a sequence of points $\crki{x_k}_{k\ge1}$ along which the left-hand side of~\eqref{Zb} diverges. We define 
    \begin{align*}
        x_k:= \sqrt{(2k\pi)^4 -1},\qquad k\geq 1.
    \end{align*}
    Hence, 
    \begin{align*}
        r(x_k) &= (2k\pi)^2 , \quad r'(x_k) = \frac{\sqrt{(2k\pi)^4 - 1}}{(2k\pi)^2 }\geq \frac{1}{\sqrt{2}},\\
        U'(r(x_k)) &= \frac{1}{2}\sin(2k\pi) + \frac{k\pi}{2}\cos(2k\pi) = \frac{k\pi}{2}.
    \end{align*}
    Consequently, 
    \begin{align*}
        \left|u'(x_k)\right|^2 =  \brk{U'(r(x_k)) r'(x_k)}^2 \geq \frac{(k\pi)^2}{8}.
    \end{align*}

    Thus, for any $\zeta,\psi\ge0$, \eqref{Zb} is violated at $x_k$ whenever $k^2\pi^2>8(\zeta^2+\psi^2)$.
    \end{proof}

\subsection{Preliminary Analysis}\label{subsec:prelim}
This subsection aims to derive an approximate descent property for $\E_t[f(x_{t})-f^\star]$, from which all convergence results in Section~\ref{sec:main} follow. The main difficulty in deriving such a property is that, within round $t$, the aggregated direction $G_t$ is affected by two sources of randomness: the stochastic samples $\xi_{i,t}^{\ell}$, $i\in\mathcal S_t$, $\ell=0,\ldots,Q-1$, used in the local updates of the selected clients, which induce stochastic gradient noise, and the random client set $\mathcal S_t$, which induces client-sampling error. These two sources are coupled, since $\mathcal S_t$ determines which samples enter $G_t$. Moreover, neither the noise nor the
client-sampling error admits a uniform bound in our setting. Under the ABC condition, the noise variance grows with the local objective gaps $f_i(x_{i,t}^{\ell})-f_i^\star$. The client-sampling error grows with the dissimilarity of the client gradients at the local iterates, which must be controlled without a BGD-type assumption.

We address these difficulties with two key ingredients. First, the virtual trajectories in Section~\ref{subsec:algorithm} separate these
two sources: conditionally on $\mathcal F_t^Q$, $G_t$ is a sample average without replacement from $n$
potential client directions, while the oracle noise is handled through $\crki{\mathcal F_t^{\ell}}_{\ell=0}^{Q}$.
Second, Lemma~\ref{AL1} controls the gradient dissimilarity through $\sigma_f^\star$, which removes the need for a BGD-type assumption.

The analysis proceeds as follows. Lemma~\ref{L2} identifies the conditional mean of $G_t$, and
Lemma~\ref{L1} quantifies the client-sampling error, namely the mean-squared deviation of the sampled average of the accumulated local gradients $\sum_{\ell=0}^{Q-1}\nabla f_i(\widetilde x_{i,t}^{\ell})$ from their average over all $n$ clients. Combining them with smoothness, Lemma~\ref{L3} derives a one-step descent inequality. Its error terms involve three quantities: the local objective gaps, the local drift, and the trajectory gradient dissimilarity. Lemmas~\ref{L5}--\ref{L4} bound these three quantities, respectively. Lemma~\ref{L7} then combines these bounds into an approximate descent property for $\mathbb E_t[f(x_{t+1})-f^\star]$, which yields the finite-time and almost sure convergence results in Section~\ref{sec:main}.

\begin{lemma}
\label{L2}
Under Assumptions \ref{A1}, \ref{A2} and \ref{A3}, the conditional mean of
\(G_t\) satisfies
\begin{equation}
    \mathbb{E}_t [G_t] = \frac{1}{n} \sum_{i=1}^n \sum_{\ell=0}^{Q-1} \mathbb{E}_t \left[ \nabla f_i(\widetilde{x}_{i,t}^\ell) \right] = Q \mathbb{E}_t\left[\nabla \bar{f}_t\right],
\end{equation}
where $\nabla \bar{f}_t := \frac{1}{nQ}\sum_{i=1}^n \sum_{\ell=0}^{Q-1} \nabla f_i(\widetilde{x}_{i,t}^\ell)$ is the virtual average gradient, i.e., the average of the local gradients along the virtual trajectories over all clients and local steps.
\end{lemma}

	\begin{proof}
		See Appendix~\ref{app:L2}.
	\end{proof}


Thus, client sampling preserves the conditional mean of the full-participation direction. However, $\mathbb E_t[G_t]$ generally differs from $Q\nabla f(x_t)$, because the local stochastic gradients are evaluated along the virtual trajectories rather than at $x_t$. This bias stems from local updates rather than client sampling, and it is controlled by the local drift (Lemma~\ref{L6}). Lemma~\ref{L1} quantifies the sampling error of the accumulated local gradients.

\begin{lemma}
\label{L1}
Suppose Assumptions \ref{A1}, \ref{A2} and \ref{A3} hold. Define the trajectory gradient dissimilarity, i.e., the gradient dissimilarity accumulated along the virtual trajectories, as
\[
    \Psi_t := \frac{1}{n}\sum_{i=1}^n \mathbb{E}_t \brk{\norm{ \sum_{\ell=0}^{Q-1} \nabla f_i(\widetilde{x}_{i,t}^\ell) - Q\nabla \bar{f}_t }^2}.
\]
Then
\begin{equation}
    \mathbb E_t\brk{
	\left\|
	\frac{1}{S}
	\sum_{i\in \mathcal{S}_t}
	\sum_{\ell=0}^{Q-1}
	\nabla f_i(\widetilde x_{i,t}^{\ell})
	- {Q\nabla \bar{f}_t}
	\right\|^2}
	=
	\frac{I_S}{S}\Psi_t.
\end{equation}
\end{lemma}

	\begin{proof}
		See Appendix~\ref{app:L3}.
	\end{proof}
Thus, the client-sampling error is proportional to $I_S/S$ and vanishes under full participation ($S=n$).

Combining Lemmas~\ref{L2} and~\ref{L1} with smoothness and the ABC condition~\eqref{ABC} yields the following one-step descent inequality.

\begin{lemma}
\label{L3}
Let Assumptions~\ref{A1}, \ref{A2} and \ref{A3} hold. Set the stepsize to satisfy $\alpha_t \leq 1/(2LQ)$ for any $t\geq 0$. Then, the iterates $\{x_t\}$ generated by Algorithm \ref{Al1} satisfy:
\begin{equation}
\label{eq:unified_recursion}
\begin{aligned}
    \mathbb{E}_t &[f(x_{t+1}) - f^\star] \leq  f(x_t) - f^\star - \frac{\alpha_t Q}{2} \|\nabla f(x_t)\|^2  \\
    &+ \frac{\alpha_t^2 C L}{Sn} \sum_{\ell=0}^{Q-1} \sum_{i=1}^n \mathbb{E}_t\left[f_i(\widetilde{x}_{i,t}^\ell) - f_i^\star\right] + \frac{\alpha_t^2 L I_S}{S}\Psi_t\\
    &+ \frac{\alpha_t L^2}{2n} \sum_{i=1}^n \sum_{\ell=0}^{Q-1} \mathbb{E}_t \brk{\norm{x_t - \widetilde{x}_{i,t}^\ell}^2}+ \frac{\alpha_t^2 L Q D}{S}.
\end{aligned}
\end{equation}
\end{lemma}

	\begin{proof}
		See Appendix~\ref{app:L4}.
	\end{proof}

Lemma~\ref{L3} reduces the convergence analysis to controlling three error terms: the sum of local objective gaps, the accumulated local drift, and the trajectory gradient dissimilarity \(\Psi_t\). We bound these quantities in Lemmas \ref{L5}--\ref{L4}, respectively.

\begin{lemma}\label{L5}
Let Assumptions \ref{A1} and \ref{A2} hold. We have
\begin{equation}\label{25}
\begin{aligned}
\sum_{i=1}^n &\sum_{\ell=0}^{Q-1} \left[ f_i(\widetilde{x}_{i,t}^\ell) - f_i^\star \right]
\leq 2nQ\left[ f(x_t) - f^\star \right] + 2nQ\sigma_f^\star + L\sum_{i=1}^n \sum_{\ell=0}^{Q-1} \left\| \widetilde{x}_{i,t}^\ell - x_t \right\|^2.
\end{aligned} 
\end{equation}
\end{lemma}

	\begin{proof}
		See \cite[Lemma 5]{huang2024distributed}. 

	\end{proof}

\begin{lemma}\label{L6}
Let Assumptions \ref{A1} and \ref{A2} hold. Set the stepsize
\[
\alpha_t \leq \min\crk{\frac{1}{2\sqrt{3}QL}, \frac{1}{2\sqrt{2QCL}} }, \quad  t\geq0.
\]
We have
\begin{equation}
\begin{aligned}\label{26}
    \sum_{\ell=0}^{Q-1} &\sum_{i=1}^n \mathbb{E}_t\left[ \left\| x_t - \widetilde{x}_{i,t}^\ell \right\|^2 \right]
    \leq 6\alpha_t^2 Q^3 n \left\| \nabla f(x_t) \right\|^2 + 4\alpha_t^2 Q^3 n(3L+2C)\sigma_f^\star + 2\alpha_t^2 Q^2 nD\\
    &+ 4\alpha_t^2 Q^3 n(2C+3L)\left[ f(x_t) - f^\star \right] .
\end{aligned} 
\end{equation}
\end{lemma}

	\begin{proof}
		See Appendix~\ref{app:L5}.
	\end{proof}

\begin{lemma}\label{L4}
Let $\Psi_t$ be the trajectory gradient dissimilarity defined in Lemma~\ref{L1}. Under Assumptions~\ref{A1}, \ref{A2} and \ref{A3}, $\Psi_t$ satisfies
\begin{equation}\label{21}
	\begin{aligned}
		    \Psi_t &\leq \frac{2 Q L^2}{n} \sum_{i=1}^n \sum_{\ell=0}^{Q-1} \mathbb{E}_t\brk{\norm{\widetilde{x}_{i,t}^\ell - x_t}^2} + 4 Q^2 L \sigma^*_f + 4 Q^2 L \brk{f(x_t) - f^\star}.
	\end{aligned}
\end{equation}
\end{lemma}

	\begin{proof}
		See Appendix~\ref{app:L7}.
	\end{proof}

Substituting the estimates from Lemmas~\ref{L5}--\ref{L4} into the recursion of Lemma~\ref{L3} yields the following approximate descent property of $\E_t[f(x_{t+1}) - f^\star]$.

\begin{lemma}
\label{L7}
Let Assumptions~\ref{A1}, \ref{A2} and~\ref{A3} hold. Set 
\begin{align*}
    \alpha_t \leq \frac{1}{8Q(L+C)},\qquad t\geq0.
\end{align*} 
Then the iterates generated by Algorithm~\ref{Al1} satisfy
\begin{equation}
    \label{eq:descent}
    \begin{aligned}
        &\E_t\brk{f(x_{t+1})-f^\star} \leq \left[1 + \frac{2\alpha_t^2 Q CL}{S} + \frac{4\alpha_t^2 Q^2 L^2 I_S}{S} + 12\alpha_t^3 Q^3 L^2(L+C)\right]\brk{f(x_t) - f^\star}
        - \frac{\alpha_t Q}{4}\norm{\nabla f(x_t)}^2 \\
        & + \frac{\alpha_t^2 LQD(1 + 2\alpha_t Q LS)}{S}+ \frac{2\alpha_t^2 Q^2 L\sigfn}{S}\prt{2L I_S + \frac{C}{Q}} + 12\alpha_t^3 Q^3L^2(L+C)\sigfn.
    \end{aligned}
\end{equation}
\end{lemma}

	\begin{proof}
		See Appendix~\ref{app:L8}.
	\end{proof}

\section{Main Results}\label{sec:main}

In this section, we establish two complementary stationarity guarantees for FedAvg under partial participation.
Theorem~\ref{T1} provides a finite-time bound on the smallest expected squared gradient norm over a prescribed horizon when a constant stepsize is used, whereas Theorem~\ref{T2} establishes convergence of the gradient norm along the full sequence, both in expectation and almost surely, under a diminishing stepsize sequence.

\begin{theorem}
\label{T1}
Let Assumptions~\ref{A1}, \ref{A2} and~\ref{A3} hold. Set \(\alpha_t=\alpha\) to satisfy
\begin{equation} \label{stepsize}
	\begin{aligned}
		        \alpha \le \min \left\{\frac{1}{8Q(L+C)}, \sqrt{\frac{S}{4TQ^2 L(2L I_S + C/Q)}},\prt{\frac{1}{24Q^3 L^2(L+C)T}}^{1/3}\right\}.
	\end{aligned}
\end{equation}

Define $\Delta_f:=f(x_0)-f^\star$. Then the iterates generated by Algorithm~\ref{Al1} satisfy
\begin{equation}\label{main}
		\begin{aligned}
			&\min_{t=0,1,\dots,T-1} \E\brk{\norm{\nabla f(x_t)}^2} \leq \frac{12\Delta_f}{\alpha Q T} + \frac{4\alpha L (D+2C\sigfn)}{S}\\
			&\qquad+ \frac{16\alpha Q L^2 I_S \sigfn}{S}+ 8\alpha^2Q^2 L^2\brk{6(L+C)\sigfn + \frac{D}{Q}}.
		\end{aligned}
	\end{equation}
In particular, if we set $\Delta_f>0$ and 
\begin{equation}
	\begin{aligned}    \label{alpha_1}
		\alpha &= \left(\sqrt{\frac{LQT( D+2C\sigma_f^\star + 4LQI_S\sigma_f^\star)}{3S\Delta_f}} + \prt{\frac{Q^3 L^2T[6(L+C)\sigfn + D/Q]}{3\Delta_f}}^{1/3} +\beta\right)^{-1}
	\end{aligned}
\end{equation}
where
\begin{equation}
	\begin{aligned}\label{co}
		\beta &= 8Q(L+C) + \sqrt{\frac{4TQ^2 L (2LI_S + C/Q)}{S}} + [24Q^3 L^2(L+C)T]^{1/3}.
	\end{aligned}
\end{equation}
Then, 
\begin{equation}
    \label{eq:minE_order}
    \begin{aligned}
        &\min_{t=0,1,\dots,T-1} \E\brk{\norm{\nabla f(x_t)}^2} =\order{\sqrt{\frac{L\Delta_f(D+C\sigfn + C\Delta_f)}{SQT}}  + \frac{(L+C)\Delta_f}{T} + L\prt{\Delta_f + \sqrt{\Delta_f\sigfn}}\sqrt{\frac{I_S}{ST}} \right.\\
        &\left.  + \prt{\frac{(L\Delta_f)^2 \Delta_f(L+C)}{T^2}}^{1/3} + \prt{\frac{(L\Delta_f)^2[(L+C)\sigfn + D/Q]}{T^2}}^{1/3}},
    \end{aligned}
\end{equation}
where $\orderi{\cdot}$ hides numerical constants.

\end{theorem}

\begin{proof}[Proof of Theorem~\ref{T1}]
    See Appendix~\ref{app:T1}.
\end{proof}

\begin{remark}
When $S=n$, we have $I_S=0$, and all sampling-induced terms in \eqref{eq:minE_order} vanish. The leading term then becomes $\orderi{\sqrt{L\Delta_f(D+C\sigma_f^\star+C\Delta_f)/(nQT)}}$, which exhibits a linear speedup in both the number of clients $n$ and the number of local updates $Q$, and has the same $\orderi{1/\sqrt{nQT}}$ dependence as the full-participation FedAvg result under the ABC condition in~\cite{huang2024distributed}.
\end{remark}

\begin{remark}
For $S<n$, both leading terms in \eqref{eq:minE_order} scale as $1/\sqrt{ST}$, so the bound retains a linear speedup in $S$. The speedup in $Q$, however, persists only when the noise-induced term dominates. Specifically, the sampling-induced term $\orderi{L\prti{\Delta_f+\sqrt{\Delta_f\sigfn}}\sqrt{I_S/(ST)}}$ in \eqref{eq:minE_order} dominates the noise-induced term $\orderi{\sqrt{L\Delta_f(D+C\sigfn+C\Delta_f)/(SQT)}}$ when $I_S\geq\nu/Q$, where
\[
    \nu:=\frac{D+C\prt{\sigfn+\Delta_f}}{L\prt{\sqrt{\Delta_f}+\sqrt{\sigfn}}^2},
\]
or equivalently,
\[
    1\leq S\leq n-\frac{\nu(n-1)}{Q}.
\]
Here, $\nu$ can be viewed as a noise-to-dissimilarity ratio: with $f(x)-f^\star$ replaced by $\Delta_f$, its numerator is the average variance bound $C\brki{f(x)-f^\star+\sigfn}+D$ implied by the ABC condition, and its denominator is the dissimilarity bound $2L\brki{f(x)-f^\star+\sigfn}$ in Lemma~\ref{AL1}, up to a factor of two. This regime is empty when $\nu>Q$.
\end{remark}

\begin{remark}
Existing analyses of FedAvg with partial participation for smooth nonconvex objectives~\cite{yang2021achieving} assume bounded variance and the BGD condition~\eqref{Zb}. Besides the stochastic gradient term $\orderi{1/\sqrt{SQT}}$, their leading-order terms contain a sampling-induced term $\orderi{I_S\zeta^2\sqrt{Q/(ST)}}$ (Table~\ref{tab:convergence_comparison}). In contrast, \eqref{eq:minE_order} holds under the ABC condition without BGD-type assumptions, and its sampling-induced terms $\orderi{L\prti{\Delta_f+\sqrt{\Delta_f\sigfn}}\sqrt{I_S/(ST)}}$ carry no $\sqrt{Q}$ factor. Hence, in our bound, more local updates do not amplify the cost of client sampling. Section~\ref{sub-3} examines this numerically.
\end{remark}

We now turn from finite-horizon guarantees to the asymptotic behavior of FedAvg under diminishing stepsizes. This complements Theorem~\ref{T1}, which controls only the best iterate over a finite horizon.

\begin{theorem}
\label{T2}
Let Assumptions~\ref{A1}, \ref{A2} and~\ref{A3} hold, and $\forall t\ge0$, let the stepsizes satisfy
\[
0<\alpha_t\le\frac{1}{8Q(L+C)},\quad
\sum_{t=0}^{\infty}\alpha_t=\infty,
\quad
\sum_{t=0}^{\infty}\alpha_t^2<\infty.
\]
Then the iterates generated by Algorithm~\ref{Al1} satisfy
\[
\lim_{t\to\infty}\E\brk{\norm{\nabla f(x_t)}}=0
\quad\text{and}\quad
\lim_{t\to\infty}\norm{\nabla f(x_t)}=0
\quad\text{a.s.}
\]
\end{theorem}

\begin{proof}[Proof of Theorem~\ref{T2}]
    See  the supplementary material.
\end{proof}

\section{Numerical Experiments}
\label{sec:simulations}

This section illustrates three implications of our analysis. Section~\ref{subsec:exp1} studies how the number of sampled clients $S$ affects convergence under the ABC condition. Section~\ref{subsec:comp} constructs instances with nearly identical BGD parameters but different values of $\sigfn$, and examines how this heterogeneity interacts with partial participation. Section~\ref{sub-3} tests whether the sampling-induced error grows with the number of local updates $Q$ when $\alpha Q$ is fixed.

\subsection{Effect of the Participation Level}
\label{subsec:exp1}

We first examine how the participation level $S$ affects the convergence of FedAvg under the ABC condition. Following~\cite{huang2024distributed}, we consider $n=16$ clients with heterogeneous local objectives
\begin{equation}\label{simu1}
	\begin{aligned}
		f_1(x) &=
		\begin{cases}
			\frac{x^2}{2}, & |x| < 1, \\
			|x| - \frac{1}{2}, & \text{otherwise},
		\end{cases}\\
		f_i(x) &= \ln\prt{1 + e^{x - i + 1}},\quad i = 2, 3, \dots, 16,
	\end{aligned}
\end{equation}
where $f_1$ is the Huber function and $\ln(1+e^x)$ is the softplus function~\cite{glorot2011deep}. The stochastic gradient oracles are
\begin{equation}
	\begin{aligned}
		g_1(x;\xi) &=
		\begin{cases}
			\nabla f_1(x) + \sqrt{|x|}, & \xi = 0, \\
			\nabla f_1(x) - \sqrt{|x|}, & \xi = 1,
		\end{cases}\\
		g_i(x;\xi) &=
		\begin{cases}
			\nabla f_i(x) + \sqrt{f_i(x)}, & \xi = 0, \\
			\nabla f_i(x) - \sqrt{f_i(x)}, & \xi = 1,
		\end{cases}\quad i = 2, 3, \dots, 16,
	\end{aligned}
\end{equation}
where $\xi\in\{0,1\}$ is a Bernoulli random variable with $\mathbb{P}(\xi = 0) = \mathbb{P}(\xi = 1) = 1/2$ drawn independently across clients, local steps, and communication rounds. Since $|x|\le f_1(x)+1/2$ and $f_i^\star=0$ for all $i$, Assumptions~\ref{A1} and~\ref{A2} hold with $L=1$, $C=1$, and $D=1/2$. In contrast, the relaxed growth condition fails \cite{huang2024distributed}.

\begin{figure}[H] 
	\centering 
	\includegraphics[width=2.5in]{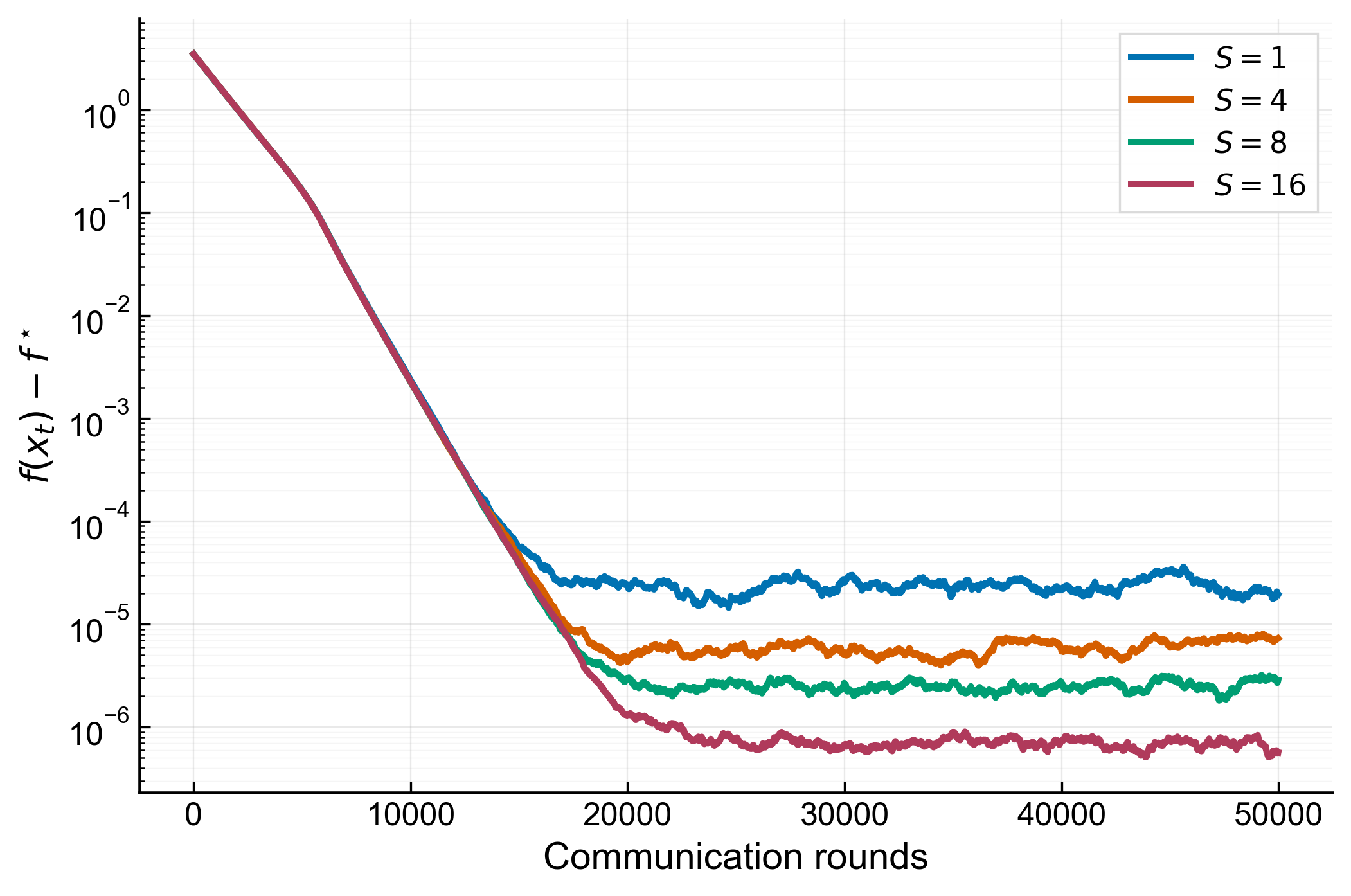} 
	\caption{Performance of FedAvg on problem~\eqref{simu1} with \(Q=5\), \(\alpha=0.001\), and different participation levels
\(S\in\{1,4,8,16\}\). The curves are averages over \(100\)
independent trials.} 
	\label{Fig.main1} 
\end{figure}

We run FedAvg with $T=50000$, \(Q=5\), $\alpha=0.001$ and \(x_0=10\), sampling $S\in\{1,4,8,16\}$ clients uniformly without replacement, where \(S=16\) corresponds to full participation. Results are averaged over \(R=100\) independent trials, and \(f^\star\) is estimated by a fine one-dimensional grid search.Figure~\ref{Fig.main1} shows convergence in terms of $f(x_t)-f^\star$ under different participation levels. 
The full-participation case \(S=16\) gives the smallest function value gap, while smaller values of \(S\) lead to larger  objective gaps. This confirms the theoretical prediction of Theorem~\ref{T1}.

\subsection{Heterogeneity Measure versus BGD Parameters}
\label{subsec:comp}

In this part, we examine whether $\sigfn$ captures changes in optimization behavior that the BGD parameters do not reflect. To this end, we construct a family of problems for which the optimized BGD sum, i.e., the smallest value of $\zeta^2+\psi^2$ such that~\eqref{Zb} holds, remains nearly constant, whereas $\sigfn$ varies substantially. We then study how this variation interacts with the participation level in FedAvg.

We consider $n=40$ clients with Huber objectives centered at random offsets $d_i$:
\begin{equation}
	\begin{aligned}
		f_i(x) &=
		\begin{cases}
			\frac{(x-d_i)^2}{2}, & |x - d_i| < 1, \\
			|x - d_i| - \frac{1}{2}, & \text{otherwise},
		\end{cases} & i &= 1,2, \dots, n,\label{simu3}
	\end{aligned}
\end{equation}
where $d_1=d_2=0$ and $d_3,\ldots,d_n$ are drawn independently from $U[0,D_{\max}]$. The parameter $D_{\max}\in\{20,50,100,200\}$ controls the minimizers $d_i$ of the local objectives. Since $f_i^\star=0$ for all $i$, we have $\sigfn=f^\star$. For each $D_{\max}$, we compute $\sigfn$ and the optimized BGD sum on the same sampled offsets.

All instances share the stochastic gradient oracle
\[
    g_i(x;\xi)
    =
    \nabla f_i(x)
    +
    \varepsilon_\xi
    \sqrt{f_i(x)-f_i^\star},
\]
where $\varepsilon_\xi$ is a random variable that satisfies $\mathbb P(\varepsilon_\xi=1)=\mathbb P(\varepsilon_\xi=-1)=1/2$. It is drawn independently across clients, local steps, communication rounds, and runs. This oracle is unbiased with variance $f_i(x)-f_i^\star$, so Assumptions~\ref{A1} and~\ref{A2} hold with $L=1$, $C=1$, and $D=0$.

We run FedAvg with $S\in\{4,20,40\}$, $Q=10$, the fixed stepsize $\alpha=0.0001$, and $T=50000$ communication rounds, where $S=40$ corresponds to full participation. All runs start from $x_0=x^\star+10$, where $x^\star$ is a minimizer of $f$, so that the initial distance to $x^\star$ is the same for all $D_{\max}$. Each configuration is repeated over $100$ independent runs.

\begin{figure*}[ht] 
	\centering 
	\includegraphics[width=0.9\textwidth]{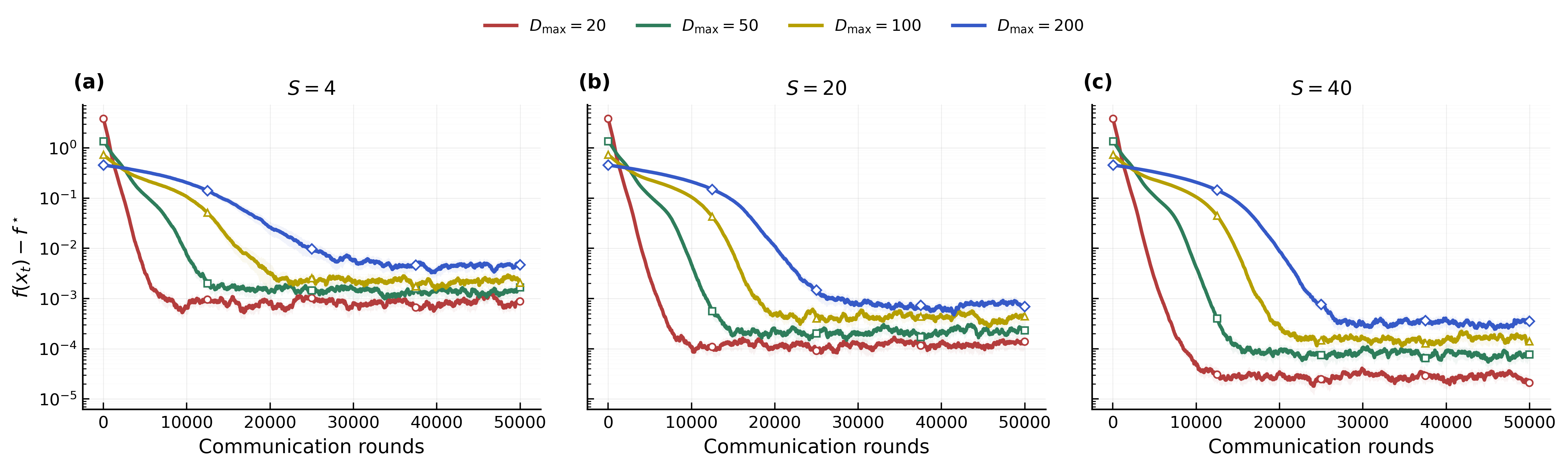} 
	\caption{Objective gap $f(x_t)-f^\star$ on problem~\eqref{simu3} for (a) $S=4$, (b) $S=20$, and (c) $S=40$, with $D_{\max}\in\{20,50,100,200\}$. } 
	\label{fig:exp2_fixedS} 
\end{figure*}

Table~\ref{tab:metrics_comparison} reports the two heterogeneity measures together with the last-iterate gap $f(x_T)-f^\star$. As $D_{\max}$ increases from $20$ to $200$, $\sigfn$ grows from $4.8$ to $52.7$, whereas the optimized BGD sum stays almost unchanged between $0.956$ and $0.998$. Meanwhile, the last-iterate gap increases with $D_{\max}$ for every $S$. Hence the change in optimization behavior is reflected by $\sigfn$ but not by the BGD parameters. 
\begin{table}[!htbp]
	\centering
	\caption{Optimized BGD parameters $\zeta^2+\psi^2$, heterogeneity measure $\sigfn$, and last-iterate gap $f(x_T)-f^\star$ on problem~\eqref{simu3} ($n=40$, $Q=10$)}
	\label{tab:metrics_comparison}
	\setlength{\tabcolsep}{2.5pt}
	\renewcommand{\arraystretch}{1.12}
	\begin{tabular*}{\columnwidth}{@{\extracolsep{\fill}}cccccc@{}}
		\toprule
		$D_{\max}$ & $\zeta^2+\psi^2$ & $\sigma_f^\star$
		& \multicolumn{3}{c}{$f(x_T)-f^\star$}\\
		\cmidrule(lr){4-6}
		& & & $S=4$ & $S=20$ & $S=40$\\
		\midrule
		20  & 0.956 & 4.838  & $9.31\times10^{-4}$ & $1.28\times10^{-4}$ & $2.76\times10^{-5}$\\
		50  & 0.998 & 12.809 & $1.34\times10^{-3}$ & $2.34\times10^{-4}$ & $7.23\times10^{-5}$\\
		100 & 0.998 & 26.096 & $2.38\times10^{-3}$ & $3.75\times10^{-4}$ & $1.65\times10^{-4}$\\
		200 & 0.998 & 52.675 & $4.51\times10^{-3}$ & $7.96\times10^{-4}$ & $3.12\times10^{-4}$\\
		\bottomrule
	\end{tabular*}
\end{table}
Figure~\ref{fig:exp2_fixedS} plots the corresponding trajectories in Table~\ref{tab:metrics_comparison}. In each subplot, a larger $D_{\max}$ slows the decrease of the objective gap and leads to a higher residual. Increasing $S$ from $4$ to $40$ lowers the residual for every $D_{\max}$. Both observations are consistent with Theorem~\ref{T1}: $\sigfn$ enters the error terms, and the sampling-induced terms decrease as $S$ increases.

\subsection{Effect of the Number of Local Updates}
\label{sub-3}

This part examines whether the local updates amplify the sampling-related error in the leading term. 
Motivated by the $Q$-dependence discussed after Theorem~\ref{T1}, we consider Problem~\eqref{simu3} with $D_{\max}=200$ and compare partial participation ($S=4$) with full participation ($S=40$). For both participation levels, we set $Q\in\{3,5,8,10,15,20,30,40,50,60,70,80,90,100,150,200\}$. For a fair comparison, the stepsize is set as $\alpha = \gamma/Q$ with a fixed constant $\gamma=10^{-2}$. Each configuration runs for $T=100000$ communication rounds and repeated over $100$ independent trials.

\begin{figure}[!htbp]
    \centering
    \includegraphics[width=0.8\textwidth]
    {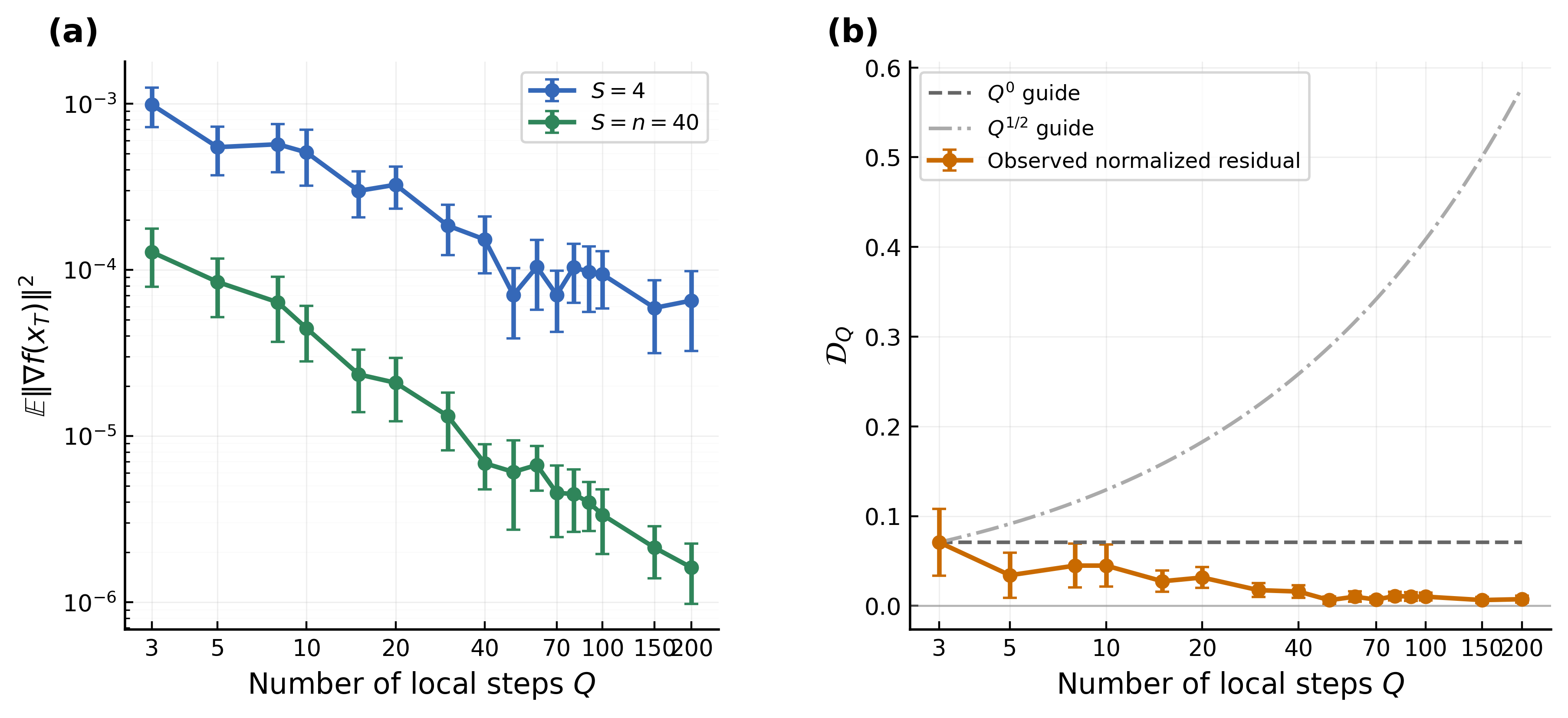}
    \caption{Dependence on the number of local updates $Q$ under $\alpha Q=\gamma$ on problem~\eqref{simu3} with $D_{\max}=200$. (a) Terminal residual $R_S(Q)=\E\brki{\normi{\nabla f(x_T)}^2}$ for $S=4$ and $S=n=40$. (b) Normalized residual $\mathcal D_Q$ with the $Q^0$ and $Q^{1/2}$ reference trends. Error bars denote $95\%$ confidence intervals over $100$ independent runs.}
    \label{fig:q_controlled}
\end{figure}

Figure~\ref{fig:q_controlled}(a) reports the terminal stationarity residual $R_S(Q) := \E\brki{\|\nabla f(x_T^{(S,Q)})\|^2}$ for $S\in\crki{4,40}$, where the expectation is estimated by averaging over the $100$ independent trials. It is observed that the residual $R_S(Q)$ decreases with increasing $Q$ for both participation levels. 
Moreover, partial participation ($S=4$) consistently yields a larger terminal residual than full participation ($S=n=40$), which is consistent with Theorem~\ref{T1}.

To further test whether the sampling-induced leading term grows as $\sqrt Q$ or is independent of $Q$, we plot the scaled residual
\begin{equation}
    \mathcal D_Q
    :=
    \frac{
        \bigl(
        \sqrt{S}R_S(Q)
        -
        \sqrt{n}R_n(Q)
        \bigr)\sqrt{T}
    }{
        L(\Delta_f +\sqrt{\Delta_f \sigma_f^\star})\sqrt{I_S}
    },
\end{equation}
against $Q$ in Figure~\ref{fig:q_controlled}(b). The difference $\sqrt{S} R_S(Q)-\sqrt{n} R_n(Q)$ is designed to cancel the leading term $\orderi{1/\sqrt{SQT}}$ in \eqref{eq:minE_order}. Then, if the sampling term is $\orderi{L(\Delta_f + \sqrt{\Delta_f\sigfn})\sqrt{I_S/(ST)}}$ as in Theorem~\ref{T1}, $\mathcal D_Q$ should be approximately constant in $Q$.
As shown in Figure~\ref{fig:q_controlled}(b), the scaled residual $\mathcal D_Q$ does not follow the $\sqrt{Q}$ trend. 
This indicates that partial participation FedAvg does not exhibit an additional $\sqrt{Q}$ amplification in the leading-order sampling term, which is consistent with the theoretical prediction of Theorem~\ref{T1}.

\section{Conclusion}
\label{sec:conclusion}

In this work, we study FedAvg with partial participation for smooth nonconvex objective functions under the ABC variance condition. Our analysis accommodates multiple local updates and uniform client sampling without replacement, without
imposing additional data heterogeneity conditions, such as the bounded gradient dissimilarity assumption. An example further demonstrates that this assumption is not implied by smoothness and lower boundedness. For a constant stepsize, we establish a finite-time bound on the smallest expected squared gradient norm over a prescribed horizon. This bound shows that partial participation introduces an additional sampling-induced error, which vanishes under full participation and, unlike existing BGD-based bounds, does not grow with the number of local updates. Under a diminishing stepsize sequence, we further establish that the gradient norm converges to zero along the full sequence, both in expectation and almost surely. Numerical simulations illustrate our theoretical findings.

\section*{Acknowledgments}

\appendix[Proofs]

Throughout the appendix, we write $g_{i,t}^{\ell}:=g_i(\widetilde x_{i,t}^{\ell};\xi_{i,t}^{\ell})$ for the stochastic gradient along the virtual trajectory of client $i$ at local step $\ell$ of round $t$, where $i\in[n]$, $t\ge0$, and $\ell\in\{0,\ldots,Q-1\}$. We also define the stochastic gradient noise $\delta_{i,t}^{\ell}:=g_{i,t}^{\ell}-\nabla f_i(\widetilde x_{i,t}^{\ell})$. By Assumption~\ref{A2}, $\condEi{\delta_{i,t}^{\ell}}{\mathcal F_t^{\ell}}=0$ and $\condEi{\normi{\delta_{i,t}^{\ell}}^2}{\mathcal F_t^{\ell}}\le C\brki{f_i(\widetilde x_{i,t}^{\ell})-f_i^\star}+D$ almost surely. In this notation, $G_t=\frac{1}{S}\sum_{i\in\mathcal S_t}\sum_{\ell=0}^{Q-1}g_{i,t}^{\ell}$.

\subsection{Proof of Lemma \ref{L2}}
\label{app:L2}

\begin{proof}
	By the definition of the update rule, $G_t$ can be expressed using the indicator variable $\mathbb{I}_{i,t}$ as $G_t = \frac{1}{S} \sum_{i=1}^n \mathbb{I}_{i,t} \sum_{\ell=0}^{Q-1} g_{i,t}^{\ell}$.  
	Since $\mathcal S_t$ is drawn uniformly from all subsets of $[n]$ with cardinality $S$, we have $\mathbb{E}_t [\mathbb{I}_{i,t}] = S/n$ for each client $i$; see \eqref{eq:indep}.
	Taking the expectation $\mathbb{E}_t$, we obtain:
	\begin{equation}
			\begin{aligned}\label{eq:EGt}
			\mathbb{E}_t [G_t] &= \mathbb{E}_t \left[ \frac{1}{S} \sum_{i=1}^n \mathbb{I}_{i,t} \sum_{\ell=0}^{Q-1} \delta_{i,t}^{\ell} \right] + \mathbb{E}_t \left[ \frac{1}{S} \sum_{i=1}^n \mathbb{I}_{i,t} \sum_{\ell=0}^{Q-1} \nabla f_i(\widetilde{x}_{i,t}^\ell) \right].
		\end{aligned}
	\end{equation}
	The first term vanishes. Indeed, the noise $\delta_{i,t}^\ell$ is $\mathcal F_t^{\ell+1}\subseteq\mathcal F_t^Q$-measurable, so the first identity in \eqref{eq:indep} removes the sampling indicator. The unbiasedness in Assumption~\ref{A2} then gives
	\begin{align*}
		\mathbb{E}_t& \left[ \frac{1}{S} \sum_{i=1}^n \mathbb{I}_{i,t} \sum_{\ell=0}^{Q-1} \delta_{i,t}^{\ell} \right]
		=\frac{1}{S}\cdot\frac{S}{n}\sum_{i=1}^n\sum_{\ell=0}^{Q-1}\mathbb E_t\left[\delta_{i,t}^{\ell}\right]\\
		&=\frac{1}{n}\sum_{i=1}^n\sum_{\ell=0}^{Q-1}\mathbb E_t\left[\condE{\delta_{i,t}^{\ell}}{\mathcal F_t^\ell}\right]
		= 0.
	\end{align*}
	Regarding the second term in \eqref{eq:EGt}, each $\nabla f_i(\widetilde{x}_{i,t}^\ell)$ is $\mathcal F_t^Q$-measurable. Applying the first identity in \eqref{eq:indep}, we obtain
	\begin{align*}
		\mathbb{E}_t &\left[ \frac{1}{S} \sum_{i=1}^n \mathbb{I}_{i,t} \sum_{\ell=0}^{Q-1} \nabla f_i(\widetilde{x}_{i,t}^\ell) \right] = \frac{1}{S} \cdot\frac{S}{n} \sum_{i=1}^n  \sum_{\ell=0}^{Q-1} \mathbb{E}_t [\nabla f_i(\widetilde{x}_{i,t}^\ell)] =Q \mathbb{E}_t\left[\nabla \bar{f}_t\right].
	\end{align*}
	Combining the identities for the two terms in \eqref{eq:EGt}, we obtain the desired result.
\end{proof}

\subsection{Proof of Lemma \ref{L1}}
\label{app:L3}

\begin{proof}
	For $i\in[n]$, define
	\begin{align*}
			b_{i,t}&:=\sum_{\ell=0}^{Q-1}\nabla f_i(\widetilde x_{i,t}^{\ell})-\frac{1}{n}\sum_{j=1}^n\sum_{\ell=0}^{Q-1}\nabla f_j(\widetilde x_{j,t}^{\ell})=\sum_{\ell=0}^{Q-1}\nabla f_i(\widetilde x_{i,t}^{\ell})-Q\nabla\bar f_t,
	\end{align*}
	so that $\Psi_t=\frac{1}{n}\sumn\E_t\brk{\norm{b_{i,t}}^2}$.
	Since \(\sumn \bbI_{i,t}=\sum_{i\in\cS_t}1\), we have
	\[
	\frac{1}{S}\sum_{i\in\cS_t}\sum_{\ell=0}^{Q-1}\nabla f_i(\widetilde x_{i,t}^{\ell})-\frac{1}{n}\sumn\sum_{\ell=0}^{Q-1}\nabla f_i(\widetilde x_{i,t}^{\ell})
	=\frac{1}{S}\sumn\bbI_{i,t}b_{i,t}.
	\]
	Hence, expanding the squared norm gives
	\[
	\begin{aligned}
		&\E_t\brk{\norm{\frac{1}{S}\sum_{i\in\cS_t}\sum_{\ell=0}^{Q-1}\nabla f_i(\widetilde x_{i,t}^{\ell})-\frac{1}{n}\sumn\sum_{\ell=0}^{Q-1}\nabla f_i(\widetilde x_{i,t}^{\ell})}^2}\\
		=\;&\frac{1}{S^2}\E_t\brk{\sumn\bbI_{i,t}\norm{b_{i,t}}^2}
		+\frac{1}{S^2}\E_t\brk{\sum_{i\ne j}\bbI_{i,t}\bbI_{j,t}\inpro{b_{i,t},b_{j,t}}}\\
		=\;&\frac{1}{S^2}\frac{S}{n}\E_t\brk{\sumn\norm{b_{i,t}}^2}
		+\frac{1}{S^2}\frac{S(S-1)}{n(n-1)}\E_t\brk{\sum_{i\ne j}\inpro{b_{i,t},b_{j,t}}}.
	\end{aligned}
	\]
	Here, the second equality uses $\bbI_{i,t}^2=\bbI_{i,t}$ and the first two identities in \eqref{eq:indep}, which apply because every $b_{i,t}$ is $\cF_t^{Q}$-measurable.
	For the cross summation, observe that $\sumn b_{i,t}=0$. Therefore,
	\[
	0=\norm{\sumn b_{i,t}}^2=\sumn\norm{b_{i,t}}^2+\sum_{i\ne j}\inpro{b_{i,t},b_{j,t}}.
	\]
	Thus,
	\[
	\begin{aligned}
		&\E_t\brk{\norm{\frac{1}{S}\sum_{i\in\cS_t}\sum_{\ell=0}^{Q-1}\nabla f_i(\widetilde x_{i,t}^{\ell})-\frac{1}{n}\sumn\sum_{\ell=0}^{Q-1}\nabla f_i(\widetilde x_{i,t}^{\ell})}^2}\\
		&=\frac{1}{S^2}\prt{\frac{S}{n}-\frac{S(S-1)}{n(n-1)}}\E_t\brk{\sumn\norm{b_{i,t}}^2}\\
		&=\frac{n-S}{S(n-1)}\E_t\brk{\frac{1}{n}\sumn\norm{b_{i,t}}^2}=\frac{I_S}{S}\Psi_t,
	\end{aligned}
	\]
	which completes the proof.
	
\end{proof}

\subsection{Proof of Lemma \ref{L3}}
\label{app:L4}

\begin{proof}
	By the $L$-smoothness of $f$ and the update rule $x_{t+1} = x_t - \alpha_t G_t$, we have:
	\begin{equation}
		\mathbb{E}_t[f(x_{t+1})] \leq f(x_t) - \alpha_t \mathbb{E}_t\brk{\inpro{\nabla f(x_t), G_t}} + \frac{L \alpha_t^2}{2} \mathbb{E}_t \brk{\norm{G_t}^2}. \label{eq:smoothness_base_new}
	\end{equation}
	
	Applying  Lemma~\ref{L2} to the inner product leads to
	\begin{equation}
		\label{eq:inner}
		\begin{aligned}
			&-\alpha_t \mathbb{E}_t\left[\inpro{\nabla f(x_t), G_t} \right]=- \alpha_t Q \mathbb{E}_t \left[ \inpro{\nabla f(x_t), \nabla \bar{f}_t} \right]\\
			&= \frac{\alpha_t Q}{2} \mathbb{E}_t \left[ \norm{\nabla f(x_t) - \nabla \bar{f}_t}^2\right] - \frac{\alpha_t Q}{2} \norm{\nabla f(x_t)}^2 - \frac{\alpha_t Q}{2}  \mathbb{E}_t \left[\norm{\nabla \bar{f}_t}^2\right]\\
			&\leq \frac{\alpha_t L^2}{2n}\sumn \sum_{\ell=0}^{Q-1}\E_t\brk{\norm{x_t - \widetilde{x}_{i,t}^{\ell}}^2} - \frac{\alpha_t Q}{2} \norm{\nabla f(x_t)}^2  - \frac{\alpha_t Q}{2}  \mathbb{E}_t \left[\norm{\nabla \bar{f}_t}^2\right],
		\end{aligned}
	\end{equation}
	where we invoked the $L$-smoothness of each $f_i$.
	
	Regarding the last term in \eqref{eq:smoothness_base_new}, it follows from the same decomposition as in \eqref{eq:EGt} that  
	\begin{equation}
		\begin{aligned}
			\mathbb{E}_t&\left[\norm{G_t}^2\right] 
			\le 2\mathbb E_t \left[
			\left\|\frac{1}{S}\sum_{i\in \mathcal{S}_t}\sum_{\ell=0}^{Q-1}\nabla f_i(\widetilde x_{i,t}^{\ell})\right\|^2\right]+2\mathbb E_t\left[\left\|\frac{1}{S}\sum_{i=1}^n\bbI_{i,t}\sum_{\ell=0}^{Q-1}
			\delta_{i,t}^{\ell}\right\|^2\right].
		\end{aligned}\label{mid}
	\end{equation}
	
	Note that $Q\nabla\bar f_t$ is $\mathcal F_t^{Q}$-measurable. By the last identity in \eqref{eq:indep},
	\begin{align*}
		\condE{\frac{1}{S}\sum_{i\in\mathcal S_t}\sum_{\ell=0}^{Q-1}\nabla f_i(\widetilde x_{i,t}^{\ell})}{\mathcal F_t^{Q}}=Q\nabla\bar f_t.
	\end{align*}
	Therefore, 
	\begin{align*}
		\E_t\brk{\inpro{Q\nabla \bar{f}_t, \frac{1}{S}\sum_{i\in \mathcal{S}_t}\sum_{\ell=0}^{Q-1}\nabla f_i(\widetilde x_{i,t}^{\ell}) - Q\nabla \bar{f}_t}} = 0.
	\end{align*}
	Then, the first term in \eqref{mid} can be bounded as follows.
	
	\begin{equation}
		\label{eq:nablaf_S}
		\begin{aligned}
			2\mathbb E_t\left[\left\|\frac{1}{S}\sum_{i\in \mathcal{S}_t}\sum_{\ell=0}^{Q-1}\nabla f_i(\widetilde x_{i,t}^{\ell})\right\|^2\right] &=2 Q^2 \mathbb E_t \left[\left\| \nabla \bar{f}_t\right\|^2\right] + 2 \mathbb E_t \left[\left\|\frac{1}{S}\sum_{i\in \mathcal{S}_t}\sum_{\ell=0}^{Q-1}\nabla f_i( \widetilde x_{i,t}^{\ell})-Q \nabla \bar{f}_t\right\|^2\right]\\
			&= 2 Q^2 \mathbb E_t \left[\left\| \nabla \bar{f}_t\right\|^2 \right]+  \frac{2 I_S }{S}\Psi_t,
		\end{aligned}
	\end{equation}
	where the second equality follows from Lemma~\ref{L1}.

	As for the second term in \eqref{mid}, each $\delta_{i,t}^{\ell}$ is $\mathcal F_t^{\ell+1}\subseteq\mathcal F_t^Q$-measurable. We first remove the sampling indicators and then show that all cross terms vanish. By the first two identities in \eqref{eq:indep}, for all $i,j\in[n]$ and $\ell,r\in\{0,\ldots,Q-1\}$,
	\begin{align*}
		\E_t\brk{\mathbb I_{i,t}\mathbb I_{j,t}\inpro{\delta_{i,t}^{\ell},\delta_{j,t}^{r}}}
		&=c_{ij} \E_t\brk{\inpro{\delta_{i,t}^{\ell},\delta_{j,t}^{r}}},\\
			c_{ij}:&=
		\begin{cases}
			\frac{S}{n}, & i=j,\\
			\frac{S(S-1)}{n(n-1)}, & i\neq j.
		\end{cases}
	\end{align*}
	If \(r>\ell\), then \(\delta_{i,t}^{\ell}\) is \(\mathcal F_t^r\)-measurable, and Assumption~\ref{A2} gives
	\[
	\E_t\brk{\inpro{\delta_{i,t}^{\ell},\delta_{j,t}^{r}}} = \E_t\brk{\inpro{\delta_{i,t}^\ell, \condE{\delta_{j,t}^r}{\mathcal{F}_t^r}}}=0.
	\]
	If \(r=\ell\) and \(i\ne j\), then \(\delta_{i,t}^{\ell}\) and \(\delta_{j,t}^{\ell}\) are conditionally independent given \(\mathcal F_t^\ell\) due to Assumption~\ref{A2}. Then,
	\[
	\E_t\brk{\inpro{\delta_{i,t}^{\ell},\delta_{j,t}^{\ell}}} = \E_t\brk{\inpro{\condE{\delta_{i,t}^\ell}{\cF_t^\ell},\condE{\delta_{j,t}^\ell}{\cF_t^\ell}}} = 0.
	\]
	Therefore, only the diagonal terms $i=j$, $r=\ell$ remain, and the ABC condition \eqref{ABC} yields
	\begin{equation}
		\label{eq:delta_var}
		\begin{aligned}
			2\E_t&\brk{\norm{\frac{1}{S}\sum_{i\in\mathcal S_t}\sum_{\ell=0}^{Q-1}\delta_{i,t}^{\ell}}^2}
			=
			\frac{2}{S^2}\cdot\frac{S}{n}\sumn\sum_{\ell=0}^{Q-1}\E_t\brk{\norm{\delta_{i,t}^{\ell}}^2}\\
			&=\frac{2}{nS}\sumn\sum_{\ell=0}^{Q-1}\E_t\brk{\condE{\norm{\delta_{i,t}^{\ell}}^2}{\cF_t^\ell}}\\
			&\leq
			\frac{2}{nS}\sumn\sum_{\ell=0}^{Q-1}
			\brk{C\E_t\brk{f_i(\widetilde x_{i,t}^{\ell})-f_i^\star}+D
			}.
		\end{aligned}
	\end{equation}

	Substituting \eqref{eq:inner}--\eqref{eq:delta_var} into \eqref{eq:smoothness_base_new} yields
	\[
	\begin{aligned}
		&\mathbb E_t[f(x_{t+1})]
		\le\;
		f(x_t)
		-
		\frac{\alpha_t Q}{2}\norm{\nabla f(x_t)}^2-
		\frac{\alpha_t Q}{2}
		\mathbb E_t \left[
		\left\|
		\nabla \bar{f}_t
		\right\|^2 \right]\\
		&+
		\frac{\alpha_t L^2}{2n}
		\sum_{i=1}^n
		\sum_{\ell=0}^{Q-1}
		\mathbb E_t \left[
		\norm{\widetilde x_{i,t}^{\ell}-x_t}^2\right]
		+
		L Q^2\alpha_t^2
		\mathbb E_t
		\left[\left\|
		\nabla \bar{f}_t
		\right\|^2\right]
			\\
		&+\frac{\alpha_t^2LI_S}{S}\Psi_t
		+
		\frac{\alpha_t^2LC}{Sn}
		\sum_{i=1}^n
		\sum_{\ell=0}^{Q-1}
		\mathbb E_t
		\brk{
			f_i(\widetilde x_{i,t}^{\ell})-f_i^\star
		}
		+
		\frac{\alpha_t^2LQD}{S}.
	\end{aligned}
	\]
	
	Letting $\alpha_t \leq 1/(2LQ)$ yields the desired result.
	
\end{proof}


\subsection{Proof of Lemma \ref{L6}}
\label{app:L5}

\begin{proof}
	Recall that the virtual local trajectory is initialized at
	$\widetilde{x}_{i,t}^{0}=x_t$ and satisfies
	\[
	\widetilde{x}_{i,t}^{\ell+1}=\widetilde{x}_{i,t}^{\ell}-\alpha_t {g_{i,t}^{\ell}},
	\qquad
	\ell=0,\ldots,Q-1.
	\]
	Unrolling the recursion gives
	\begin{equation}
		\begin{aligned}
			\widetilde{x}_{i,t}^{\ell} - x_t
			= -\alpha_t \sum_{j=0}^{\ell-1}\brk{\delta_{i,t}^{j} + \nabla f_i(\widetilde{x}_{i,t}^{j}) }.
		\end{aligned}
	\end{equation}
	Hence,
	\begin{equation}
		\label{eq:xitl_diff}
		\begin{aligned}
			\frac{1}{2\alpha_t^2}\E_t\brk{\norm{\widetilde{x}_{i,t}^{\ell}-x_t}^2}&\leq
			\E_t\brk{\norm{\sum_{j=0}^{\ell-1}\delta_{i,t}^{j}}^2}
			+\E_t\brk{\norm{\sum_{j=0}^{\ell-1}\nabla f_i(\widetilde{x}_{i,t}^{j})}^2}\\
			&\leq \E_t\brk{\norm{\sum_{j=0}^{\ell-1}\delta_{i,t}^{j}}^2} + \ell \sum_{j=0}^{\ell-1}\E_t\brk{\norm{\nabla f_i(\widetilde{x}_{i,t}^{j})}^2}.
		\end{aligned}
	\end{equation}
	Following similar derivations as in \eqref{eq:delta_var}, the first term in \eqref{eq:xitl_diff} can be bounded as follows:
	\begin{equation}
		\begin{aligned}
			\E_t\brk{\norm{\sum_{j=0}^{\ell-1}\delta_{i,t}^{j}}^2}
			&=\sum_{j=0}^{\ell-1}\E_t\brk{\norm{\delta_{i,t}^{j}}^2}\leq C\sum_{j=0}^{\ell-1}\E_t\brk{f_i(\widetilde{x}_{i,t}^{j})-f_i^\star}+\ell D.
		\end{aligned}
	\end{equation}
	Hence,
	\begin{equation}
		\label{eq:xitl_diff_s1}
		\begin{aligned}
			\frac{1}{2\alpha_t^2}&\E_t\brk{\norm{\widetilde{x}_{i,t}^{\ell}-x_t}^2}
			\leq\;
			C\sum_{j=0}^{\ell-1}\E_t\brk{f_i\prti{\widetilde{x}_{i,t}^{j}}-f_i^\star}+\ell D
			+\ell\sum_{j=0}^{\ell-1}\E_t\brk{\norm{\nabla f_i\prti{\widetilde{x}_{i,t}^{j}}}^2}.
		\end{aligned}
	\end{equation}
	
	Summing \eqref{eq:xitl_diff_s1} over $i=1,\ldots,n$ and
	$\ell=0,\ldots,Q-1$ yields
	\begin{equation}
		\label{eq:xitl_diff_S2}
		\begin{aligned}
			&\sum_{\ell=0}^{Q-1}\sumn\E_t\brk{\norm{\widetilde{x}_{i,t}^{\ell}-x_t}^2} 
			\leq \alpha_t^2Q^2\sum_{\ell=0}^{Q-1}\sumn\E_t\brk{\norm{\nabla f_i\prti{\widetilde{x}_{i,t}^{\ell}}}^2} + \alpha_t^2Q^2nD+ 2\alpha_t^2QC\sum_{\ell=0}^{Q-1}\sumn\E_t\brk{f_i\prti{\widetilde{x}_{i,t}^{\ell}}-f_i^\star}\\
			&\leq \alpha_t^2 Q^2 nD + 4\alpha_t^2 Q^2 n C\sigfn + 4\alpha_t^2 Q^2n C\brk{f(x_t) - f^\star} + 2\alpha_t^2 QCL\sumn \sum_{\ell=0}^{Q-1}\E_t\brk{\norm{\widetilde{x}_{i,t}^{\ell}-x_t}^2} \\
			&\quad+ \alpha_t^2Q^2\sum_{\ell=0}^{Q-1}\sumn\E_t\brk{\norm{\nabla f_i\prti{\widetilde{x}_{i,t}^{\ell}}}^2},
		\end{aligned}
	\end{equation}
	where we invoked Lemma~\ref{L5}.

	We now bound the last term in \eqref{eq:xitl_diff_S2}. Noting 
	\begin{equation}
		\begin{aligned}
			\nabla f_i&\prti{\widetilde{x}_{i,t}^{\ell}}
			=
			\brk{\nabla f_i\prti{\widetilde{x}_{i,t}^{\ell}}-\nabla f_i(x_t)}
			+\brk{\nabla f_i(x_t)-\nabla f(x_t)}
			+\nabla f(x_t),
		\end{aligned}
	\end{equation}
	it follows from the $L$-smoothness of $f_i$ that 
	\begin{equation}
		\label{eq:nabla_fi_sum}
		\begin{aligned}
			&\sum_{\ell=0}^{Q-1}\sumn\E_t\brk{\norm{\nabla f_i\prti{\widetilde{x}_{i,t}^{\ell}}}^2}
			\leq 3L^2\sum_{\ell=0}^{Q-1}\sumn\E_t\brk{\norm{\widetilde{x}_{i,t}^{\ell}-x_t}^2} +3Qn\normi{\nabla f(x_t)}^2+3Q\sumn\norm{\nabla f_i(x_t)-\nabla f(x_t)}^2\\
			&\quad\leq 3L^2\sum_{\ell=0}^{Q-1}\sumn\E_t\brk{\norm{\widetilde{x}_{i,t}^{\ell}-x_t}^2}
			+3Qn\normi{\nabla f(x_t)}^2 + 6QnL \brk{f(x_t) - f^\star} + 6QnL \sigfn,
		\end{aligned}
	\end{equation}
	where we invoked Lemma~\ref{AL1}.
	
	Substituting \eqref{eq:nabla_fi_sum} into \eqref{eq:xitl_diff_S2} and rearranging the terms, we obtain
	\begin{equation}
		\begin{aligned}
			&\prt{1-3\alpha_t^2Q^2L^2-2\alpha_t^2QCL}
			\sum_{\ell=0}^{Q-1}\sumn\E_t\brk{\norm{\widetilde{x}_{i,t}^{\ell}-x_t}^2}\\
			&\leq
			3\alpha_t^2Q^3n\normi{\nabla f(x_t)}^2
			+2\alpha_t^2Q^3n(2C+3L)\brk{f(x_t)-f^\star}
			+2\alpha_t^2Q^3n(3L+2C)\sigfn
			+\alpha_t^2Q^2nD.
		\end{aligned}
	\end{equation}
	
	Invoking the stepsize condition yields the desired result.
\end{proof}

\subsection{Proof of Lemma \ref{L4}}
\label{app:L7}

\begin{proof}
	For any vectors $\crki{a_i}_{i=1}^n$ with mean $\bar a=\frac{1}{n}\sumn a_i$ and any $c\in\R^d$, we have $\frac{1}{n}\sumn\normi{a_i-\bar a}^2\le\frac{1}{n}\sumn\normi{a_i-c}^2$. Applying this with $a_i=\sum_{\ell=0}^{Q-1}\nabla f_i(\widetilde x_{i,t}^{\ell})$, for which $\bar a=Q\nabla\bar f_t$, and $c=Q\nabla f(x_t)$, and then taking $\E_t$, gives
	\begin{equation}\label{eq:Psi_s1}
		\begin{aligned}
			\Psi_t&=\frac{1}{n}\sumn\E_t\brk{\norm{\sum_{\ell=0}^{Q-1}\nabla f_i(\widetilde x_{i,t}^{\ell})-Q\nabla\bar f_t}^2}\le\frac{1}{n}\sumn\E_t\brk{\norm{\sum_{\ell=0}^{Q-1}\nabla f_i(\widetilde x_{i,t}^{\ell})-Q\nabla f(x_t)}^2}.
		\end{aligned}
	\end{equation}
	Noting
	\begin{align*}
	\sum_{\ell=0}^{Q-1}&\nabla f_i(\widetilde x_{i,t}^{\ell})-Q\nabla f(x_t)
	=\sum_{\ell=0}^{Q-1}\brk{\nabla f_i(\widetilde x_{i,t}^{\ell})-\nabla f_i(x_t)}+Q\brk{\nabla f_i(x_t)-\nabla f(x_t)},
	\end{align*}
	we obtain
	\[
	\begin{aligned}
		\Psi_t&\le\frac{2}{n}\sumn\E_t\brk{\norm{\sum_{\ell=0}^{Q-1}\brk{\nabla f_i(\widetilde x_{i,t}^{\ell})-\nabla f_i(x_t)}}^2}+\frac{2Q^2}{n}\sumn\norm{\nabla f_i(x_t)-\nabla f(x_t)}^2\\
		&\leq \frac{2QL^2}{n}\sumn\sum_{\ell=0}^{Q-1}\E_t\brk{\norm{\widetilde x_{i,t}^{\ell}-x_t}^2} + 4Q^2L\brk{f(x_t)-f^\star} + 4Q^2L\sigfn,
	\end{aligned}
	\]
	where we invoked the Cauchy--Schwarz inequality and the $L$-smoothness of $f_i$ for the first term, and Lemma~\ref{AL1} for the second term. This completes the proof.
\end{proof}

\subsection{Proof of Lemma \ref{L7}}
\label{app:L8}
\begin{proof}
	
	Substituting the estimates from Lemmas~\ref{L5} and \ref{L4} into \eqref{eq:unified_recursion} yields 
	\begin{equation}
		\label{eq:F_descent_s1}
		\begin{aligned}
			\E_t\brk{f(x_{t+1})-f^\star}
			&\leq \brk{1 + \frac{2\alpha_t^2 Q CL}{S} + \frac{4\alpha_t^2 Q^2 L^2 I_S}{S}}\brk{f(x_t) - f^\star} \\
			& + \frac{\alpha_t L^2}{2n}\prt{1 + \frac{2\alpha_t C}{S} + \frac{4\alpha_t L Q I_S}{S}}\sumn\sum_{\ell=0}^{Q-1}\E_t\brk{\norm{\widetilde x_{i,t}^{\ell}-x_t}^2}\\
			& + \frac{\alpha_t^2 LQD}{S} + \frac{2\alpha_t^2 Q^2 L\sigfn}{S}\prt{2L I_S + \frac{C}{Q}}- \frac{\alpha_t Q}{2}\norm{\nabla f(x_t)}^2.
		\end{aligned}
	\end{equation}
	Noting that $\alpha_t \leq 1/[8Q(L+C)]$ and $I_S \leq 1$, it follows that 
	\begin{align*}
		1 + \frac{2\alpha_t C}{S} + \frac{4\alpha_t Q L I_S}{S}\leq 2,\qquad t\geq 0.
	\end{align*}
	Invoking \eqref{26} leads to 
	\begin{equation}
		\label{eq:F_descent_s2}
		\begin{aligned}
			&\E_t\brk{f(x_{t+1})-f^\star} \leq \left[1 + \frac{2\alpha_t^2 Q CL}{S} + 12\alpha_t^3 Q^3 L^2(L+C) + \frac{4\alpha_t^2 Q^2 L^2 I_S}{S}\right]\brk{f(x_t) - f^\star}  + \frac{\alpha_t^2 LQD(1 + 2\alpha_t Q LS)}{S}\\
			& + \frac{2\alpha_t^2 Q^2 L\sigfn}{S}\prt{2L I_S + \frac{C}{Q}} + 12\alpha_t^3 Q^3L^2(L+C)\sigfn - \frac{\alpha_t Q}{2}\prt{1 - 12\alpha_t^2 Q^2 L^2}\norm{\nabla f(x_t)}^2.
		\end{aligned}
	\end{equation}

	Letting $\alpha_t\leq 1/(2\sqrt{6}QL)$ yields the desired result.
	
\end{proof}

\subsection{Proof of Theorem \ref{T1}}
\label{app:T1}

\begin{proof}
	
	The proof of Theorem \ref{T1} relies on Lemma \ref{L7} and the following result from \cite{mishchenko2020random}.
	
	\begin{lemma}[Lemma 6, \cite{mishchenko2020random}]\label{L9}
		Suppose that there exist constants $a,b,c \geq 0$ and nonnegative sequences $\{s_t\}_{t=0}^T$, $\{q_t\}_{t=0}^T$ such that for any $t$ satisfying $0 \leq t \leq T$, we have the recursion
		$$
		s_{t+1} \leq (1+a)s_t - bq_t + c.
		$$
		Then, the following holds:
		\[
		\min_{t=0,1,\dots,T-1} q_t \leq \frac{(1+a)^T}{bT} s_0 + \frac{c}{b}.
		\]
	\end{lemma}
	
	To apply Lemma \ref{L9}, we set
	\begin{align*}
		s_t &= \E[f(x_t) - f^\star],\quad q_t = \E\brk{\norm{\nabla f(x_t)}^2},\quad b = \frac{\alpha Q}{4}\\
		c &= \frac{\alpha^2 LQD(1 + 2\alpha Q LS)}{S}+ \frac{2\alpha^2 Q^2 L\sigfn}{S}\prt{2L I_S + \frac{C}{Q}} + 12\alpha^3 Q^3L^2(L+C)\sigfn,\\
		a &= \frac{2\alpha^2 Q CL}{S} + \frac{4\alpha^2 Q^2 L^2 I_S}{S} + 12\alpha^3 Q^3 L^2(L+C).
	\end{align*}
	
	Letting 
	\begin{align*}
		\alpha \leq &\min\left\{\sqrt{\frac{S}{4TQ^2 L(2L I_S + C/Q)}}, \prt{\frac{1}{24Q^3 L^2(L+C)T}}^{1/3} \right\}
	\end{align*}
	leads to 
	\begin{align*}
		\prt{1 + a}^T\leq \exp(aT)\leq \exp(1)\leq 3.
	\end{align*}
	Therefore, from Lemma \ref{L9}, we have
	\begin{equation}
		\label{eq:minE_s1}
		\begin{aligned}
			&\min_{t=0,1,\dots,T-1} \E\brk{\norm{\nabla f(x_t)}^2} \leq \frac{12\Delta_f}{\alpha Q T} + \frac{4\alpha L D}{S}\\
			&+ \frac{8\alpha Q L\sigfn}{S}\prt{2L I_S + \frac{C}{Q}}+ 8\alpha^2Q^2 L^2\brk{6(L+C)\sigfn + \frac{D}{Q}},
		\end{aligned}
	\end{equation}
	which proves \eqref{main}. 
	
	For the choice of $\alpha$ in \eqref{alpha_1}, we have 
	\begin{equation}
		\label{eq:alpha_inv}
		\begin{aligned}
			\frac{\Delta_f}{\alpha QT} &= \sqrt{\frac{L\Delta_f(D+2C\sigfn + 4LQI_S\sigma_f^\star)}{3SQT}}  + \sqrt{\frac{4L\Delta_f(2L\Delta_f I_S + C\Delta_f /Q)}{ST}} \\
			&+ \prt{\frac{24(L\Delta_f)^2 \Delta_f(L+C)}{T^2}}^{1/3} + \frac{8(L+C)\Delta_f}{T}+ \prt{\frac{(L\Delta_f)^2[6(L+C)\sigfn + D/Q]}{3T^2}}^{1/3},\\
			\alpha^2 Q^2 L^2 &\leq \prt{\frac{3L\Delta_f}{T[6(L+C)\sigfn + D/Q]}}^{2/3}.
		\end{aligned}
	\end{equation}
	Substituting \eqref{eq:alpha_inv} into \eqref{eq:minE_s1} yields \eqref{eq:minE_order}.

\end{proof}

\subsection{Proof of Theorem \ref{T2}}
\label{app:T2}
\begin{proof}
	The proof relies on Lemma~\ref{L7} and the following result, which we apply with $\Phi:=\nabla f$, $x^k:=x_t$, $\mu_k:=\alpha_t$, and the filtration $\crki{\cF_t}_{t\ge0}$.
	
	\begin{lemma}[{\cite[Theorem~2.1]{li2022unified}}]\label{L10}
		Let the mapping $\Phi:\R^d\to\R^m$ and the sequences $\crki{x^k}_{k\ge0}\subseteq\R^d$ and $\crki{\mu_k}_{k\ge0}\subseteq\R_{++}$ be given. Consider the following conditions.
		\begin{description}
			\item[(P.1)] The function $\Phi$ is $L_\Phi$-Lipschitz continuous for some $L_\Phi>0$, i.e., $\normi{\Phi(x)-\Phi(y)}\le L_\Phi\normi{x-y}$ for all $x,y\in\R^d$.
			\item[(P.2)] There exists a constant $a>0$ such that $\sum_{k=0}^\infty\mu_k\E\brki{\normi{\Phi(x^k)}^a}<\infty$.
		\end{description}
		The following statements are valid.
		\begin{enumerate}
			\item[i)] Let (P.1)--(P.2) hold and suppose further that
			\begin{description}
				\item[(P.3)] there exist constants $A,B,b\ge0$ and $p_1,p_2,q>0$ such that $\E\brki{\normi{x^{k+1}-x^k}^q}\le A\mu_k^{p_1}+B\mu_k^{p_2}\E\brki{\normi{\Phi(x^k)}^b}$;
				\item[(P.4)] $\crki{\mu_k}_{k\ge0}$ is bounded, $\sum_{k=0}^\infty\mu_k=\infty$, $a,q\ge1$, $a\ge b$, and $p_1,p_2\ge q$.
			\end{description}
			Then $\lim_{k\to\infty}\E\brki{\normi{\Phi(x^k)}}=0$.
			\item[ii)] Let (P.1)--(P.2) hold and assume further that
			\begin{description}
				\item[(P.3')] there exist constants $A,b\ge0$, $p_1,p_2,q>0$ and random vectors $\boldsymbol{A}_k,\boldsymbol{B}_k$ such that $x^{k+1}=x^k+\mu_k^{p_1}\boldsymbol{A}_k+\mu_k^{p_2}\boldsymbol{B}_k$, where, for all $k$, $\boldsymbol{A}_k$ and $\boldsymbol{B}_k$ are $\cF_{k+1}$-measurable, $\condE{\boldsymbol{A}_k}{\cF_k}=0$ almost surely, $\E\brki{\normi{\boldsymbol{A}_k}^q}\le A$, and $\limsup_{k\to\infty}\normi{\boldsymbol{B}_k}^q/\prti{1+\normi{\Phi(x^k)}^b}<\infty$ almost surely;
				\item[(P.4')] $\mu_k\to0$, $\sum_{k=0}^\infty\mu_k=\infty$, $\sum_{k=0}^\infty\mu_k^{2p_1}<\infty$, $q\ge2$, $qa\ge b$, $p_1>\frac{1}{2}$, and $p_2\ge1$.
			\end{description}
			Then $\lim_{k\to\infty}\normi{\Phi(x^k)}=0$ almost surely.
		\end{enumerate}
	\end{lemma}
	
	\paragraph{Verification of (P.1)} Since each $f_i$ is $L$-smooth, so is $f$. Hence, $\Phi=\nabla f$ is $L$-Lipschitz continuous.
	
	\paragraph{Verification of (P.2)} Since $\alpha_t\le 1/[8Q(L+C)]$, we have $12\alpha_t^3Q^3L^2(L+C)\le\frac{3}{2}\alpha_t^2Q^2L^2$ and $2\alpha_tQLS\le S/4$. Hence, \eqref{eq:descent} in Lemma~\ref{L7} gives
	\begin{equation}\label{eq:T2_rec}
		\begin{aligned}
    	\E_t&\brk{f(x_{t+1})-f^\star}\le\prt{1+K_1\alpha_t^2}\brk{f(x_t)-f^\star}-\frac{\alpha_tQ}{4}\norm{\nabla f(x_t)}^2+K_2\alpha_t^2,
		\end{aligned}
	\end{equation}
	where
\begin{align*}
	K_1&:=\frac{2QL}{S}\prt{C+2QLI_S}+\frac{3}{2}Q^2L^2,\\
	K_2&:=LQD\prt{\frac{1}{S}+\frac{1}{4}}+K_1\sigfn
\end{align*}
	are independent of $t$. Since $\sum_{t=0}^\infty\alpha_t^2<\infty$, the Robbins--Siegmund theorem~\cite{robbins1971convergence} implies that $f(x_t)$ converges almost surely to a finite random variable. Taking total expectations in \eqref{eq:T2_rec} and discarding the gradient term gives $\E\brki{f(x_{t+1})-f^\star}\le\prti{1+K_1\alpha_t^2}\E\brki{f(x_t)-f^\star}+K_2\alpha_t^2$. Iterating this inequality and using $1+u\le \exp(u)$ yields
	\begin{align*}
	&\sup_{t\ge0}\E\brk{f(x_t)-f^\star}\le M:=\exp\prt{K_1\sum_{t=0}^\infty\alpha_t^2}\prt{\Delta_f+K_2\sum_{t=0}^\infty\alpha_t^2}<\infty.
	\end{align*}
	Taking total expectations in \eqref{eq:T2_rec} again, summing over $t=0,\ldots,T$, and using $\E\brki{f(x_{T+1})-f^\star}\ge0$, we obtain
	\[
	\frac{Q}{4}\sum_{t=0}^{T}\alpha_t\E\brk{\norm{\nabla f(x_t)}^2}\le\Delta_f+\prt{K_1M+K_2}\sum_{t=0}^\infty\alpha_t^2.
	\]
	Letting $T\to\infty$ shows that (P.2) holds with $a=2$.
	
	By \eqref{mid}, \eqref{eq:nablaf_S} and \eqref{eq:delta_var},
	\begin{align*}
			\E_t&\brk{\norm{G_t}^2}\le2Q^2\E_t\brk{\norm{\nabla\bar f_t}^2}+\frac{2I_S}{S}\Psi_t+\frac{2}{nS}\sumn\sum_{\ell=0}^{Q-1}\prt{C\E_t\brk{f_i(\widetilde x_{i,t}^{\ell})-f_i^\star}+D}.
	\end{align*}
	Adding and subtracting $\nabla f_i(x_t)$ inside $\nabla\bar f_t$, and using Jensen's inequality and the $L$-smoothness of $f_i$, we have 
	\begin{align*}
		2Q^2&\E_t\brk{\norm{\nabla\bar f_t}^2}\le4Q^2\norm{\nabla f(x_t)}^2+\frac{4QL^2}{n}\sumn\sum_{\ell=0}^{Q-1}\E_t\brk{\norm{\widetilde x_{i,t}^{\ell}-x_t}^2}.
	\end{align*} 
	Substituting these estimates, together with Lemma~\ref{L4} for $\Psi_t$ and Lemma~\ref{L5}, we obtain
	\[
	\begin{aligned}
		\E_t&\brk{\norm{G_t}^2}\le{}4Q^2\norm{\nabla f(x_t)}^2+\frac{2L}{n}\prt{2QL+\frac{2QLI_S+C}{S}}\sumn\sum_{\ell=0}^{Q-1}\E_t\brk{\norm{\widetilde x_{i,t}^{\ell}-x_t}^2}\\
		&+\frac{4Q}{S}\prt{2QLI_S+C}\brk{f(x_t)-f^\star+\sigfn}+\frac{2QD}{S}.
	\end{aligned}
	\]
	Since $I_S\le1$ and $S\ge1$, the coefficient of the drift term is at most $2L(4QL+C)/n$. Moreover, $\alpha_t\le1/(8Q(L+C))$ implies $2L(4QL+C)\alpha_t^2Q^2\le Q/8$. Hence, by Lemma~\ref{L6},
	\begin{align*}
		&\frac{2L(4QL+C)}{n}\sumn\sum_{\ell=0}^{Q-1}\E_t\brk{\norm{\widetilde x_{i,t}^{\ell}-x_t}^2}\le\frac{3Q^2}{4}\norm{\nabla f(x_t)}^2\\
		&+\frac{Q^2(2C+3L)}{2}\brk{f(x_t)-f^\star+\sigfn}+\frac{QD}{4}.
	\end{align*}
	Combining the two displays with $\normi{\nabla f(x_t)}^2\le2L\brki{f(x_t)-f^\star}$ yields
	\begin{equation}\label{eq:T2_G}
		\E_t\brk{\norm{G_t}^2}\le K_3+K_4\brk{f(x_t)-f^\star},
	\end{equation}
	where
	\[
	\begin{aligned}
		K_3&:=\brk{\frac{4Q}{S}\prt{2QLI_S+C}+\frac{Q^2}{2}(2C+3L)}\sigfn+QD\prt{\frac{2}{S}+\frac{1}{4}},\\
		K_4&:=\frac{19}{2}Q^2L+\frac{4Q}{S}\prt{2QLI_S+C}+\frac{Q^2}{2}(2C+3L)
	\end{aligned}
	\]
	are independent of $t$.
	Consequently, $\sup_{t\ge0}\E\brki{\normi{G_t}^2}\le K_3+K_4M<\infty$.
	
	\paragraph{Verification of (P.3) and (P.4)} Since $x_{t+1}-x_t=-\alpha_tG_t$, we have
	\[
	\E\brk{\norm{x_{t+1}-x_t}^2}=\alpha_t^2\E\brk{\norm{G_t}^2}\le\prt{K_3+K_4M}\alpha_t^2.
	\]
	Thus, (P.3) holds with $q=b=p_1=p_2=2$, $A=K_3+K_4M$, and $B=0$. Condition (P.4) holds because $\crki{\alpha_t}$ is bounded, $\sum_{t=0}^\infty\alpha_t=\infty$, $a=q=2$, $a\ge b$, and $p_1,p_2\ge q$. Therefore, Lemma~\ref{L10}~i) gives $\lim_{t\to\infty}\E\brki{\normi{\nabla f(x_t)}}=0$.
	
	\paragraph{Verification of (P.3') and (P.4')} Write $x_{t+1}=x_t+\alpha_tA_t+\alpha_tB_t$ with $A_t:=\E_t\brki{G_t}-G_t$ and $B_t:=-\E_t\brki{G_t}$, so that $p_1=p_2=1$. Both $A_t$ and $B_t$ are $\cF_{t+1}$-measurable, $\E_t\brki{A_t}=0$, and $\E\brki{\normi{A_t}^2}\le\E\brki{\normi{G_t}^2}\le K_3+K_4M$. Moreover, it follows from Jensen's inequality and \eqref{eq:T2_G} that 
	\begin{align*}
		\norm{B_t}^2\le\E_t\brk{\norm{G_t}^2}\le K_3+K_4\brk{f(x_t)-f^\star}.
	\end{align*}
	Since $f(x_t)$ converges almost surely to a finite limit, $\sup_{t\ge0}\normi{B_t}^2<\infty$ almost surely, and hence $\limsup_{t\to\infty}\normi{B_t}^2/\prti{1+\normi{\nabla f(x_t)}^2}<\infty$ almost surely. Thus, (P.3') holds with $q=b=2$. Finally, (P.4') holds because $\alpha_t\to0$, $\sum_{t=0}^\infty\alpha_t=\infty$, $\sum_{t=0}^\infty\alpha_t^2<\infty$, $q=2$, $qa=4\ge b$, $p_1=1>\frac{1}{2}$, and $p_2=1$. Therefore, Lemma~\ref{L10}~ii) yields $\lim_{t\to\infty}\normi{\nabla f(x_t)}=0$ almost surely.
\end{proof}


\bibliographystyle{IEEEtran}
\bibliography{reference}

\newpage

 




\vfill

\end{document}